\documentclass{amsart}

\usepackage[all]{xy}

\numberwithin{equation}{section}

\theoremstyle{plain}
\newtheorem{theorem}{Theorem}[section]
\newtheorem{corollary}[theorem]{Corollary}
\newtheorem{lemma}[theorem]{Lemma}
\newtheorem{proposition}[theorem]{Proposition}
\theoremstyle{remark}
\newtheorem{remark}[theorem]{Remark}
\theoremstyle{definition}
\newtheorem{definition}[theorem]{Definition}

\makeatletter

\def\R{\mathbb{R}}
\def\C{\mathbb{C}}
\def\Z{\mathbb{Z}}

\newcommand{\lie}[1]{\mathfrak{#1}}
\newcommand{\lier}[1]{\mathfrak{#1}_{0}}
\newcommand{\roots}[2]{\Delta({#1}, {#2})}
\newcommand{\proots}[2]{\Delta^{+}({#1}, {#2})}
\newcommand{\stWh}[2]{I_{\eta,{#1},{#2}}^{\circ}}
\newcommand{\stWhmg}[2]{I_{\eta,{#1},{#2}}}
\newcommand{\stWhL}[1]{\stWh{\Lambda}{#1}}
\newcommand{\stWhLmg}[1]{\stWhmg{\Lambda}{#1}}
\def\stWhLs{\stWhL{\sigma}}
\def\stWhLsmg{\stWhLmg{\sigma}}
\def\stWhLg{\stWhL{\gamma}}
\def\stWhLgmg{\stWhLmg{\gamma}}
\def\WhL{I_{\eta,\Lambda}^{\circ}}
\def\gK{\lie{g},K}
\def\brgK{(\gK)}
\def\rhom{\rho_{\lie{m}}}
\def\lietm{\lie{t}_{\lie{m}}}

\def\Hom{\mathrm{Hom}}

\def\tr{\mathrm{tr}}

\def\Ind{\mathrm{Ind}}
\def\Ad{\mathrm{Ad}}
\def\ad{\mathrm{ad}}
\def\pr{\mathrm{pr}}
\def\Dim{\mathrm{Dim}}

\def\Wh{\mathrm{Wh}}
\def\sgn{\mathrm{sgn}}

\newcommand{\vect}[1]{\mathbf{#1}}
\def\I{\sqrt{-1}}

\allowdisplaybreaks[3]

\title[Standard Whittaker $\brgK$-modules]
{On the composition series of the standard Whittaker
  $\brgK$-modules}
\author{Kenji Taniguchi}
\address{
Department of Physics and Mathematics, 
Aoyama Gakuin University, 
5-10-1, Fuchinobe, Chuo-ku, Sagamihara, Kanagawa 252-5258, Japan. }
\email{taniken@gem.aoyama.ac.jp}

\begin{document}

\begin{abstract}
For a real reductive linear Lie group $G$, 
the space of Whittaker functions
is the representation space induced from a non-degenerate unitary
character of the Iwasawa nilpotent subgroup. 
Defined are the standard Whittaker $(\mathfrak{g},K)$-modules, 
which are $K$-admissible submodules of the space of
Whittaker functions. 
We first determine the structures of them when the infinitesimal
characters characterizing them are generic. 
As an example of the integral case, we determine the composition
series of the standard Whittaker $\brgK$-module 
when $G$ is the group $U(n,1)$ and the infinitesimal character 
is regular integral. 
\end{abstract}

\maketitle

\section{Introduction}
\label{section:introduction}

One of the most basic problems in representation theory is to study the
composition series of a standard representation. 
In the category of highest weight modules, 
Verma modules play the role
of standard representations, 
and in the category of
Harish-Chandra modules, principal series representations do. 
The composition series problem is called the
Kazhdan-Lusztig conjecture. 

In this paper, the author propose a Whittaker version of standard
$\brgK$-modules and study their composition series problem. 

Let $G$ be a real reductive linear Lie group in the sense of \cite{V2} 
and $G=KAN$ be an Iwasawa decomposition of it. 
Let $\eta: N \longrightarrow \C^{\times}$ be a unitary character of $N$ 
and denote the differential representation $\lier{n} \to \I \R$ 
of it by the same letter $\eta$. 
We assume $\eta$ is non-degenerate, i.e. it is
non-trivial on every root space corresponding to a simple root of
$\Delta^{+}(\lier{g}, \lier{a})$. 
Define 
\begin{equation}\label{eq:space of Whittaker functions}
C^{\infty}(G/N; \eta) 
:= \{f: G \overset{C^{\infty}}{\longrightarrow} \C \,|\, 
f(gn) = \eta(n)^{-1} f(g), \enskip g \in G, n \in N\}
\end{equation}
and call it the space of Whittaker functions on $G$. 
This is a representation space of
$G$ by the left translation, which is denoted by $L$. 
Let $C^{\infty}(G/N;\eta)_{K}$ be the subspace of
$C^{\infty}(G/N;\eta)$ consisting of $K$-finite vectors. 
As for the subrepresentations of this space, 
there are many deep and interesting results, 
called the theory of Whittaker models. 
On the other hand, it is not too much to say that the structure of
the whole space is not known at all. 
Though our ultimate goal is to determine the structure of 
$C^{\infty}(G/N; \eta)$, this space is too large to analyze. 
So we need to cut off a submodule of suitable size from it. 

Let, as usual, $M$ be the centralizer of $A$ in $K$, and let 
\[
M^{\eta} 
:= 
\{m \in M |\, \eta(m^{-1}nm) = \eta(n), n \in N\}
\]
be the stabilizer of $\eta$ in $M$. 
This subgroup acts naturally on $C^{\infty}(G/N; \eta)_{K}$ by the
right translation. 
Consider the subspace of $C^{\infty}(G/N; \eta)_{K}$ consisting of
those functions $f$ which satisfy the following conditions: 
\begin{enumerate}
\item
$f$ is a joint eigenfunction of $Z(\lie{g})$ (the center of
  the universal enveloping algebra $U(\lie{g})$) 
with eigenvalue 
$\chi_{\Lambda}$: 
$L(z) f = \chi_{\Lambda}(z) f$, $z \in Z(\lie{g})$. 
\item
For an irreducible representation $(\sigma, V_{\sigma}^{M^{\eta}})$ 
of $M^{\eta}$, 
$f$ is in the $\sigma^{\ast}$-isotypic subspace 
($\sigma^{\ast}$ is the dual of $\sigma$) 
with respect to the right action of $M^{\eta}$. 
\item
$f$ grows moderately at the infinity (\cite{W}).  
\end{enumerate}
Denote by $\WhL$ the subspace consisting of 
$f \in C^{\infty}(G/N; \eta)_{K}$ satisfying (1). 
Then 
\[
\WhL
\simeq 
\bigoplus_{\sigma \in \widehat{M^{\eta}}} 
\Hom_{M^{\eta}} 
(V_{\sigma^{\ast}}^{M^{\eta}}, 
\WhL) 
\otimes V_{\sigma^{\ast}}^{M^{\eta}},
\]
and the space 
$\Hom_{M^{\eta}} 
(V_{\sigma^{\ast}}^{M^{\eta}}, 
\WhL)$ is isomorphic to 
\begin{align}
\stWhLs 
= C^{\infty}(G/M^{\eta}N; 
\sigma \otimes&\eta)_{K,\Lambda}
\notag\\
:= 
\{f: G \overset{C^{\infty}}{\longrightarrow}
V_{\sigma}^{M^{\eta}} 
|\, 
&
f(gmn) = \eta(n)^{-1} \sigma(m)^{-1} f(g), 
\enskip g \in G, m \in M^{\eta}, n \in N; 
\notag\\
& L(z) f = \chi_{\Lambda}(z) f, \enskip z \in Z(\mathfrak{g}); 
\enskip 
\mbox{left $K$-finite} 
\}.
\notag
\end{align}
Therefore, the space of functions $f$ satisfying the above conditions
(1)--(3) is isomorphic to 
\begin{align}
\stWhLsmg
:=
\{f\in \stWhLs \,|\,  f \mbox{ grows moderately at the infinity}\}. 
\notag
\end{align}
We call these the {\it standard Whittaker $\brgK$-modules}. 
Note that these are {\it not} 
the ``standard Whittaker module'' defined in \cite{Kos}. 
It is easy to show that 
these are $K$-admissible and then have finite length 
(Corollary~\ref{corollary:K-admissibility}).

Though the composition series problem of standard Whittaker
$\brgK$-modules is an interesting problem by itself, 
we may hope to apply the result of it
to the analysis of principal series representations. 
$\stWhLsmg$ is induced from $M^{\eta}N$ and the behavior of 
$f \in \stWhLsmg$ on $A$ is controlled by the infinitesimal
character and the asymptotic behavior, so we may
think this module is near to the principal series representation. 
Therefore, it is significant to compare the structure of
this module and that of a principal series representation. 
According to the theory of Whittaker models, 
an irreducible Harish-Chandra module $\pi$ can be a submodule of 
$\stWhLsmg$ only if the Gelfand-Kirillov dimension of $\pi$
is equal to $\dim N$ (\cite{M1}). 
On the other hand, any irreducible Harish-Chandra module can be a
submodule of some principal series representation. 
This difference comes from the difference of the structures of 
$\stWhLsmg$ and principal series. 
So if you understand the common features and the different points of
these modules, then new insights of Whittaker models and principal
series are expected to be obtained. 

For example, assume $G=SL(2,\R)$ and the infinitesimal character
$\Lambda$ is regular dominant integral.  
In this case 
$M = M^{\eta} \simeq \{\pm 1\}$, 
so we identify an irreducible representation of 
$M$ and that of $M^{\eta}$, which is denoted by $\sigma$. 
There are two irreducible representations of $\{\pm 1\}$, one is
trivial, denoted by $1$, and the other is the signature
representation, denoted by $-1$. 
Let 
$\rho_{A} = \frac{1}{2} \tr (\ad_{\lie{a}}|_{\mathrm{Lie}(N)}) 
\in \mathrm{Lie}(A)^{\ast}$. 
Then the principal series representation 
$\Ind_{MAN}^{G}(\sigma \otimes e^{\Lambda+\rho_{A}})$ is reducible if
and only if $\Lambda \equiv \sigma+1 \mod\, 2$. 
There are four equivalence classes of irreducible Harish-Chandra
modules with the infinitesimal character $\Lambda$. 
The irreducible principal series is denoted by
$\overline{\pi}_{01}^{-}$, 
the irreducible finite dimensional representation by
$\overline{\pi}_{01}^{+}$ 
and two discrete series are denoted by $\pi_{0}$, $\pi_{1}$. 
It is well known that the composition series of reducible principal
series are 
\begin{align*}
& 
\Ind_{MAN}^{G}(\sigma \otimes e^{\Lambda+\rho_{A}})
\quad \simeq \quad 
\begin{xy}
(-4,-4)*{\pi_{0}}="A_{1}",
(0,4)*{\overline{\pi}_{01}^{+}}="A_{2}",
(4,-4)*{\pi_{1}}="A_{3}",
\ar "A_{2}";"A_{1}"
\ar "A_{2}";"A_{3}"
\end{xy} 
& 
& 
\Ind_{MAN}^{G}(\sigma \otimes e^{-\Lambda+\rho_{A}})
\quad \simeq \quad 
\begin{xy}
(4,4)*{\pi_{1}}="A_{1}",
(0,-4)*{\overline{\pi}_{01}^{+}}="A_{2}",
(-4,4)*{\pi_{0}}="A_{3}",
\ar "A_{1}";"A_{2}"
\ar "A_{3}";"A_{2}"
\end{xy} 
\end{align*}
For the meaning of these diagrams, 
see Definition~\ref{definition:diagram of composition series}. 
On the other hand, 
if $\sigma \in \widehat{M^{\eta}} = \widehat{M}$ corresponds to
the reducible (resp. irreducible) principal series, then the
composition series of the standard Whittaker $\brgK$-modules are 
\begin{align*}
& 
\stWhLsmg
\quad \simeq \quad 
\begin{xy}
(0,7)*{\pi_{0}}="A_{1}",
(0,0)*{\overline{\pi}_{01}^{+}}="A_{2}",
(0,-7)*{\pi_{1}}="A_{3}",
\ar "A_{1}";"A_{2}"
\ar "A_{2}";"A_{3}"
\end{xy} 
& 
& 
I_{-\eta,\Lambda,\sigma}
\quad \simeq \quad 
\begin{xy}
(0,7)*{\pi_{1}}="A_{1}",
(0,0)*{\overline{\pi}_{01}^{+}}="A_{2}",
(0,-7)*{\pi_{0}}="A_{3}",
\ar "A_{1}";"A_{2}"
\ar "A_{2}";"A_{3}"
\end{xy} 
& 
&
(\mbox{resp. } 
\quad 
\stWhLsmg
\simeq 
\overline{\pi}_{01}^{-})
\end{align*}
for an appropriately chosen $\eta$. 
This result can be obtained by direct computation. 

In this paper, we first determine the structure of 
$I_{\eta, \Lambda, \sigma}$ when $\Lambda$ is generic. 
Let $X_{P}(\delta,\nu)$ be the Harish-Chandra module of the
$C^{\infty}$-induced principal series
representation 
$C^{\infty}\mbox{-} 
\Ind_{P}^{G}(\delta \otimes e^{\nu+\rho_{A}})$. 
Here $\delta$ is an irreducible representation of $M$. 
The Weyl group of $\lie{g}$ is denoted by $W$. 
Let $\mathcal{H}_{\Lambda}$ be the set of equivalence classes of
irreducible Harish-Chandra modules with the infinitesimal character
$\Lambda$. 
We call $\Lambda$ generic if every principal series representation
with the infinitesimal character $\Lambda$ is irreducible. 
The main theorem on the generic case is 
\begin{theorem}[Theorem~\ref{theorem:generic case}]
Let $G$ be a real reductive linear Lie group. 
Suppose $\Lambda$ is generic and $\sigma$ is an irreducible
representation of $M^{\eta}$. 
Then $\WhL$ is completely reducible. 
Moreover, the irreducible decomposition of $\stWhLsmg$ is
given by 
\begin{align}
&
\stWhLsmg
\simeq 
\bigoplus_{X_{P}(\delta,\nu) \in \mathcal{H}_{\Lambda}}
m_{\delta}(\sigma) 
X_{P}(\delta,\nu), 
&
&
m_{\delta}(\sigma) 
= \dim \Hom_{M^{\eta}}(\delta|_{M^{\eta}},\sigma).
\end{align}
\end{theorem}

For the non-generic case, 
there is little result that can be applied to general groups. 
Therefore, we examine the case $G=U(n,1)$ in the second half of this
paper so that it becomes a springboard to the study of general cases. 
Let  $\overline{\pi}_{i,j}$ be the irreducible
Harish-Chandra module of $U(n,1)$ defined in 
\S~\ref{subsection:classification of (g,K)-modules}. 
The main result on this case is 
\begin{theorem}[Theorem~\ref{theorem:main}]
Suppose $G = U(n,1)$ and the infinitesimal character $\Lambda$ is
regular integral. 
If the highest weight of 
$\sigma \in \widehat{M^{\eta}} 
\simeq \widehat{U(n-2)} \times \widehat{U(1)}$ satisfies  
\eqref{eq:condition for sigma} for some $i=1,\dots,n-1$,
$j=2,\dots,n+1-i$, 
then the composition series of $I_{\eta,\Lambda,\sigma}$ is 
\begin{equation*}
\stWhLsmg
\quad \simeq \quad 
\begin{xy}
(0,12)*{\overline{\pi}_{i-1,j+1}}="A_{11}",
(15,12)*{\overline{\pi}_{i-1,j-1}}="A_{12}",
(30,12)*{\overline{\pi}_{i+1,j+1}}="A_{13}",
(45,12)*{\overline{\pi}_{i+1,j-1}}="A_{14}",
(0,0)*{\overline{\pi}_{i-1,j}}="A_{21}",
(15,0)*{\overline{\pi}_{i,j+1}}="A_{22}",
(30,0)*{\overline{\pi}_{i,j-1}}="A_{23}",
(45,0)*{\overline{\pi}_{i+1,j}}="A_{24}",
(22.5,-12)*{\overline{\pi}_{i,j}}="A_{31}",
\ar "A_{11}";"A_{21}"
\ar "A_{11}";"A_{22}"
\ar "A_{12}";"A_{21}"
\ar "A_{12}";"A_{23}"
\ar "A_{13}";"A_{22}"
\ar "A_{13}";"A_{24}"
\ar "A_{14}";"A_{23}"
\ar "A_{14}";"A_{24}"
\ar "A_{21}";"A_{31}"
\ar "A_{22}";"A_{31}"
\ar "A_{23}";"A_{31}"
\ar "A_{24}";"A_{31}"
\end{xy}
\end{equation*}
Here, if $i+j=n$ or $n+1$, the modules 
$\overline{\pi}_{a,b}$, $a+b>n+1$, are regarded to be zero and
the arrows starting from or ending at such modules are omitted. 
\end{theorem}

This paper is organized as follows. 
The generic case is treated in \S\ref{section:generic case}. 
The main result of this section is 
Theorem~\ref{theorem:generic case}. 
From \S\ref{section:U(n,1)} later, we put $G=U(n,1)$ and
examine the composition series of $\stWhLsmg$ 
when the infinitesimal character is regular integral. 
\S\ref{section:U(n,1)} recalls the structure of $U(n,1)$ and the
classification of irreducible Harish-Chandra modules of it. 
In \S\ref{section:composition factors}, 
we first show that $I_{\eta,\Lambda,\sigma}$ has a unique irreducible
submodule if it is non-zero. 
Also determined are 
the possible irreducible modules appearing in the composition series
of it. 
In \S\ref{section:determination of the composition series}, 
the composition series of $\stWhLsmg$ is completely
determined. 
For this step, we use the explicit form
of $K$-type shift operators and the central elements of the
universal enveloping algebra. 
The key lemmas for our calculation are 
Lemma~\ref{lemma:condition for vanishing under shift} and 
\ref{lemma:central element coming from shift}, 
and the main theorem of the latter half of this paper is 
Theorem~\ref{theorem:main}. 
In \S\ref{section:ending remark}, another formulation of our problem
is discussed.

Before going ahead, we introduce notation used in this paper. 
For a real Lie group $L$, the Lie algebra of it is denoted by
$\lier{l}$ and its complexification by 
$\lie{l} = \lier{l} \otimes_{\R} \C$. 
This notation will be applied to groups denoted by other Roman letters
in the same way without comment. 
For a compact Lie group $L$, the set of equivalence classes of
irreducible representations of $L$ is denoted by $\widehat{L}$. 
The representation space of $\pi \in \widehat{L}$ is denoted by
$V_{\pi}^{L}$. 
When $L$ is connected and $\pi$ is the irreducible representation
whose highest weight is $\lambda$, we also denote it by 
$V_{\lambda}^{L}$. 
For $\pi \in \widehat{L}$, the contragredient representation is
denoted by $\pi^{\ast}$, and if $\lambda$ is the highest weight of
$\pi$, then the highest weight of $\pi^{\ast}$ is denoted by
$\lambda^{\ast}$. 

Suppose that $K$ is a maximal compact subgroup of a real
reductive group $G$. 
For a $\brgK$-module $\pi$, the $K$-spectrum 
$\{\tau \in \widehat{K}\,|\, \tau \subset \pi|_{K}\}$ is denoted by 
$\widehat{K}(\pi)$.

For a numerical vector 
$\vect{a}=(a_{1},\dots,a_{\ell}) \in \C^{\ell}$ or $\R^{\ell}$, 
write $\vert \vect{a} \vert := \sum_{i=1}^{\ell} a_{i}$. 
This notation will be applied for an element of the dual of a Cartan
subalgebra when this space is identified with numerical vector space
by using some fixed basis. 

The author would like to thank Hiroshi Yamashita, Kyo Nishiyama, 
Noriyuki Abe and Hisayosi Matumoto for helpful discussion on this problem. 
He also thanks T\^{o}ru Umeda, Minoru Itoh and Akihito Wachi for
useful advice on the determinant type central element of the
universal enveloping algebra. 
This research is partially supported by JSPS Grant-in-Aid Scientific
Research (C) \# 19540226.


\section{The generic case}\label{section:generic case}

In this section, we first write down the differential equations
characterizing $\stWhLs$. 
After that, we determine the structure of the standard Whittaker
$\brgK$-modules when $\Lambda$ is generic. 
This is the first main theorem of this paper. 
As a corollary to the proof of this theorem, 
the $K$-admissibility of any standard Whittaker $\brgK$-module is obtained.

The $K$-type decomposition of $\stWhLs$ is given by 
\begin{align*}
\stWhLs
&\simeq 
\bigoplus_{\tau \in\widehat{K}} 
\Hom_{K}(V_{\tau}^{K}, \stWhLs) 
\otimes V_{\tau}^{K}. 
\end{align*}
By Iwasawa decomposition, 
an element of $\Hom_{K}(V_{\tau}^{K}, \stWhLs)$ is
determined by its restriction to $A$. 
For $\phi_{1} \in \Hom_{K}(V_{\tau}^{K}, \stWhLs)$, 
$a \in A$, $m \in M^{\eta}$ and $v \in V_{\tau}^{K}$, 
\[
\phi_{1}(\tau(m)v)(a) 
= 
L(m) (\phi_{1}(v))(a) 
= 
\phi_{1}(v)(m^{-1}a) 
= \phi_{1}(v)(a m^{-1}) 
= \sigma(m) \phi_{1}(v)(a).
\]
Therefore, we may identify $\phi_{1}$ with an element
$\phi_{2}$ of 
$C^{\infty}
(A \rightarrow \Hom_{M^{\eta}}(V_{\tau}^{K},
V_{\sigma}^{M^{\eta}}))$ 
by $\phi_{1}(v)(a) = \phi_{2}(a)(v)$, 
$v \in V_{\tau}^{K}, a \in A$. 
The $\lie{g}$ action on $\phi_{1}$ can be transferred to $\phi_{2}$, 
which we denote by $X \cdot \phi_{2}$: 
$(X \cdot \phi_{2})(a)(v) = L(X)(\phi_{1}(v))(a)$. 
It follows that $\Hom_{K}(V_{\tau}^{K}, \stWhLs)$ is
isomorphic to 
\begin{align}\label{eq:tau isotypic in I}
\stWhLs(\tau) 
:= \{\phi_{2} \in 
C^{\infty}
&(A \rightarrow \Hom_{M^{\eta}}(V_{\tau}^{K}, V_{\sigma}^{M^{\eta}}))
\,|\, 
z \cdot \phi_{2} = \chi_{\Lambda}(z) \phi_{2}, 
z \in Z(\lie{g})\}, 
\end{align}
and the subspace $\stWhLsmg(\tau)$ of $\stWhLs(\tau)$ consisting of
functions which grow moderately at the infinity is isomorphic to 
$\Hom_{K}(V_{\tau}^{K}, \stWhLsmg)$.

We write down the action of $z \in Z(\lie{g})$ on $\phi_{2}$. 
Firstly, the $U(\lie{k})$ action is 
\[
(u \cdot \phi_{2})(a)(v) 
= (L(u) \phi_{1}(v))(a)
= \phi_{1}(\tau(u)v)(a) 
= \phi_{2}(a) (\tau(u)v), 
\quad 
u \in U(\lie{k}). 
\]

Secondly, consider the action of $U(\lie{a})$. 
Denote by $\Pi = \{\alpha_{1}, \dots, \alpha_{l}\}$ the set of
simple roots of $\proots{\lier{g}}{\lier{a}}$. 
Let $H_{1}, \dots, H_{l}$ be the basis of $\lier{a}$ dual to
$\alpha_{1}, \dots, \alpha_{l}$: 
$\alpha_{i}(H_{j}) = \delta_{ij}$. 
For $t = (t_{1}, \dots, t_{l}) \in (\R_{>0})^{l}$, define 
\[
a_{t} = \exp \left(-\sum_{i=1}^{l} (\log t_{i}) H_{i} \right) 
\in A.
\]
For 
$u = \sum_{\vect{n}} c_{\vect{n}} 
H_{1}^{n_{1}} \cdots H_{l}^{n_{l}} \in U(\lie{a})$, 
$\vect{n} = (n_{1}, \dots, n_{l}) \in (\Z_{\geq 0})^{l}$, 
$c_{\vect{n}} \in \C$, 
let 
\begin{equation}\label{eq:d_t}
\partial_{t}(u) 
= 
\sum_{\vect{n}} c_{\vect{n}} 
\partial_{1}^{n_{1}} \cdots \partial_{l}^{n_{l}}, 
\quad 
\partial_{i} 
:= t_{i} \frac{\partial}{\partial t_{i}}. 
\end{equation}
Then 
\[
(u \cdot \phi_{2}(a_{t}))(v) 
= (\partial_{t}(u) \phi_{2}(a_{t}))(v), 
\quad u \in U(\lie{a}). 
\]

Lastly, we write down the action of $U(\lie{n})$. 
Denote by $(\lier{g})_{\alpha}$ the root space corresponding
to a root $\alpha$. 
Let $\{N_{\alpha,j} \,|\, \alpha \in \proots{\lier{g}}{\lier{a}}, 
1 \leq j \leq \dim (\lier{g})_{\alpha}\}$ be a basis of
$\lier{n}$ such that it satisfies 
$\eta(N_{\alpha,j}) \not= 0$ if $\alpha \in \Pi$ and $j=1$, and 
$\eta(N_{\alpha,j}) = 0$ otherwise. 
We define 
\begin{equation}\label{eq:eta_t}
\eta_{t}(N_{\alpha,j}) 
= 
\begin{cases}
t_{i} \eta(N_{\alpha_{i},1}) \quad 
\mbox{if } \alpha=\alpha_{i} \in \Pi, j=1,
\\
0 \quad \mbox{otherwise}
\end{cases}
\end{equation}
and extend it to an algebra homomorphism 
$U(\lie{n}) \to \C[t_{1}, \dots, t_{l}]$. 
Then 
\[
(u \cdot \phi_{2}(a_{t}))(v) 
= (L(u) \phi_{1}(v))(a_{t}) 
= \eta_{t}(u) \phi_{1}(v)(a_{t}) 
= \eta_{t}(u) \phi_{2}(a_{t})(v), 
\]
for $u \in U(\lie{n})$. 
Therefore, $U(\lie{g})$ acts on 
$C^{\infty}
(A \rightarrow \Hom_{M^{\eta}}(V_{\tau}^{K}, V_{\sigma}^{M^{\eta}}))$ 
by 
\begin{align}\label{eq:action of Ug}
&((u_{\lie{n}} u_{\lie{a}} u_{\lie{k}}) \cdot \phi_{2})(a_{t})(v) 
= 
\eta_{t}(u_{\lie{n}}) 
\partial(u_{\lie{a}}) 
\phi_{2}(a_{t})(\tau(u_{\lie{k}}) v), 
\\
&
u_{\lie{n}} \in U(\lie{n}), 
u_{\lie{a}} \in U(\lie{a}), 
u_{\lie{k}} \in U(\lie{k}). 
\notag
\end{align}

Choose a Cartan subalgebra $\lietm$ of $\lie{m}$. 
Fix a positive system $\proots{\lie{m}}{\lietm}$ 
of the root system $\roots{\lie{m}}{\lietm}$. 
Let $\delta \in \widehat{M}$. 
Its highest weight with respect to $\proots{\lie{m}}{\lietm}$ 
is denoted by $\mu_{\delta}$. 
Note that since we assume every Cartan subgroup of $G$ is commutative 
\cite[(0.1.2) f)]{V2}, the highest weight $\mu_{\delta}$ of the
restriction of $\delta$ to the identity component of $M$ is well
defined even if $M$ is not connected. 
Let $P=MAN$ be the minimal parabolic subgroup of $G$ corresponding to
our Iwasawa $N$. 
For $\delta \in \widehat{M}$ and $\nu \in \lie{a}^{\ast}$, 
let $X_{P}(\delta,\nu)$ be the
Harish-Chandra module of the smooth principal series
representation 
$C^{\infty}\mbox{-} 
\Ind_{P}^{G}(\delta \otimes e^{\nu+\rho_{A}})$. 

\begin{definition}\label{definition:generic pair}
An infinitesimal character $\Lambda$ is 
called {\it generic} if every principal series representation
$X_{P}(\delta,\nu)$ which admits the infinitesimal character $\Lambda$
is irreducible. 
\end{definition}

Choose a Cartan subalgebra $\lie{h} := \lietm+\lie{a}$ of
$\lie{g}$. 
Let $W=W(\lie{g},\lie{h})$ and 
$W_{\lie{m}} = W(\lie{m},\lietm)$ be the Weyl groups of
$\lie{g}$ and $\lie{m}$, respectively. 
The little Weyl group is denoted by $W(G,A)$. 
It is well known that the infinitesimal character of 
$X_{P}(\delta,\nu)$ is $\Lambda \in \lie{h}^{\ast}$ if and
only if $(\mu_{\delta}+\rhom, \nu)$ is in the orbit 
$W \cdot \Lambda$. 
It is also well known that two principal series representations 
$X_{P}(\delta,\nu)$ and $X_{P}(\delta',\nu')$ have the same
composition factors if and only if there exists $w \in W(G,A)$ such
that $(\delta',\nu') = (w \cdot \delta, w \cdot \nu)$. 
We denote by $\mathcal{A}_{\Lambda}$ the set of 
$(\delta,\nu) \in \widehat{M} \times \lie{a}^{\ast}$ satisfying 
$(\mu_{\delta}+\rhom, \nu) \in W \cdot \Lambda$. 
The set of equivalence classes of irreducible Harish-Chandra modules
is denoted by $\mathcal{H}_{\Lambda}$. 
Note that, if $\Lambda$ is generic, 
every member of $\mathcal{H}_{\Lambda}$ is a principal series
representation. 
Therefore, $\mathcal{H}_{\Lambda}$ is parametrized by 
$W(G,A) \backslash \mathcal{A}_{\Lambda}$, the set of 
$W(G,A)$ -orbits in $\mathcal{A}_{\Lambda}$. 

The first main result of this paper is the following theorem. 
\begin{theorem}\label{theorem:generic case}
Suppose $\Lambda$ is generic and $\sigma$ is an irreducible
representation of $M^{\eta}$. 
Then $\WhL$ is completely reducible. 
Moreover, the irreducible decomposition of $\stWhLsmg$ is
given by 
\begin{align}
&
\stWhLsmg
\simeq 
\bigoplus_{X_{P}(\delta,\nu) \in \mathcal{H}_{\Lambda}}
m_{\delta}(\sigma) 
X_{P}(\delta,\nu), 
&
&
m_{\delta}(\sigma) 
= \dim \Hom_{M^{\eta}}(\delta|_{M^{\eta}},\sigma).
\label{eq:composition series, generic-2}
\end{align}
\end{theorem}
\begin{proof}
We first count the dimension of $\stWhLs(\tau)$. 

Let $\bar{\lie{n}}$ be the nilpotent subalgebra opposite to
$\lie{n}$. 
We denote by $\lie{u}_{\lie{m}}$ and $\bar{\lie{u}}_{\lie{m}}$ 
the nilpotent subalgebras in $\lie{m}$ corresponding to 
$\proots{\lie{m}}{\lietm}$ and
$-\proots{\lie{m}}{\lietm}$, respectively. 
Then $\lie{u} := \lie{u}_{\lie{m}} + \bar{\lie{n}}$ is the
nilradical of a Borel subalgebra $\lie{h} + \lie{u}$. 
Let $\rhom$ be half the sum of elements in 
$\proots{\lie{m}}{\lietm}$. 
We define non-shifted Harish-Chandra maps $\gamma_{1}'$, $\gamma_{2}'$
and $\gamma'$ by 
\begin{align*}
& \gamma_{1}' : 
U(\lie{g}) = U(\lie{m}+\lie{a}) \oplus 
(\lie{n} U(\lie{g})+U(\lie{g}) \bar{\lie{n}})
\rightarrow U(\lie{m}+\lie{a}) 
\\
& \gamma_{2}' : 
U(\lie{m}+\lie{a}) = U(\lie{h}) \oplus 
(\bar{\lie{u}}_{\lie{m}} U(\lie{m}+\lie{a})
+U(\lie{m}+\lie{a}) \lie{u}_{\lie{m}})
\rightarrow U(\lie{h}) 
\\
& \gamma' = \gamma_{2}' \circ \gamma_{1}' : 
U(\lie{g}) = U(\lie{h}) \oplus 
(\bar{\lie{u}} U(\lie{g})+U(\lie{g}) \lie{u})
\rightarrow U(\lie{h}), 
\end{align*}
respectively. 
Then Harish-Chandra maps are given by 
\begin{align*}
& \gamma_{1} = \tau_{1} \circ \gamma_{1}' 
: Z(\lie{g}) \rightarrow Z(\lie{m} + \lie{a}), 
\quad 
\tau_{1}(H) = H + \rho_{A}(H), H \in \lie{a}, 
\\
& \gamma_{2} = \tau_{2} \circ \gamma_{2}' 
: Z(\lie{m}+\lie{n}) \overset{\sim}{\rightarrow} U(\lie{h})^{W}, 
\quad 
\tau_{2}(H) = H -\rhom(H), H \in \lie{h}, 
\\
& \gamma = \gamma_{2} \circ \gamma_{1} 
: Z(\lie{g}) \overset{\sim}{\rightarrow} U(\lie{h})^{W},
\end{align*}
respectively. 
The infinitesimal character $\chi_{\Lambda}$ is, of course, 
defined by 
$\chi_{\Lambda}(z) = \gamma(z)(\Lambda)$, $z \in Z(\lie{g})$. 

Choose a $K$-type $\tau \in \widehat{K}(\stWhLs)$. 
Suppose 
$z = \sum_{p} u_{\lie{n}}^{(p)} u_{\lie{a}}^{(p)} u_{\lie{k}}^{(p)} 
\in U(\lie{n}) \otimes U(\lie{a}) \otimes U(\lie{k})$ is an element of
$Z(\lie{g})$. 
By \eqref{eq:action of Ug}, 
elements $\phi_{2}$ of $\stWhLs(\tau)$
are characterized by the system of differential equations 
\begin{equation}\label{eq:diff eq for phi_2}
\sum_{p} 
\eta_{t}(u_{\lie{n}}^{(p)}) 
\partial(u_{\lie{a}}^{(p)}) 
\phi_{2}(a_{t})(\tau(u_{\lie{k}}^{(p)}) v)
= \chi_{\Lambda}(z) \phi_{2}(a_{t})(v), 
\quad a_{t} \in A, v \in V_{\tau}^{K}. 
\end{equation}
We know that the system of partial differential equations 
$z \cdot \phi_{2} = \chi_{\Lambda}(z) \phi_{2}$, $z \in Z(\lie{g})$, 
has regular singularity at $t=0$ 
(see the definitions \eqref{eq:d_t}, \eqref{eq:eta_t} of 
$\partial_{t}(u)$ and $\eta_{t}$). 
Suppose there exists a non-zero solution $\phi_{2}$ of this system. 
Then its leading term $\phi_{2}^{0}$ satisfies the system of
differential equations 
\begin{equation}\label{eq:equation for leading term}
\sum_{p} 
\eta_{0}(u_{\lie{n}}^{(p)}) 
\partial(u_{\lie{a}}^{(p)}) 
\phi_{2}^{0}(a_{t})(\tau(u_{\lie{k}}^{(p)}) v)
= \chi_{\Lambda}(z) \phi_{2}^{0}(a_{t})(v), 
\ 
z = \sum_{p} u_{\lie{n}}^{(p)} u_{\lie{a}}^{(p)}
u_{\lie{k}}^{(p)} 
\in Z(\lie{g}). 
\end{equation}
Since 
$\sum_{p} 
\eta_{0}(u_{\lie{n}}^{(p)}) u_{\lie{a}}^{(p)} u_{\lie{k}}^{(p)}$ is
equivalent to $z$ modulo $\lie{n} U(\lie{g})$, 
it is the same as $\gamma_{1}'(z)$. 
Therefore, we may assume that, 
if $\eta_{0}(u_{\lie{n}}^{(p)}) \not= 0$, 
then $u_{\lie{k}}^{(p)} \in
Z(\lie{m})$. 
Suppose $\delta \in \widehat{M}$ satisfies 
$\sigma \subset \delta|_{M^{\eta}}$ and 
$\delta \subset \tau|_{M}$. 
Choose a non-zero vector $v$ in the
$\delta$-isotypic subspace of $V_{\tau}^{K}$. 
Then $u_{\lie{k}}^{(p)} v = 
(\mu_{\delta}+\rhom)(\gamma_{2}(u_{\lie{k}}^{(p)})) v$.  
Let $\nu+\rho_{A} \in \lie{a}^{\ast} \simeq \C^{l}$ be a
characteristic exponent of a solution $\phi_{2}$ to 
\eqref{eq:diff eq for phi_2}. 
Here, we identified $\lie{a}^{\ast}$ with $\C^{l}$ by 
$\sum_{i=1}^{l} c_{i} \alpha_{i} \leftrightarrow 
(c_{1},\dots,c_{l})$. 
The equation \eqref{eq:equation for leading term} says that $\mu_{\delta}$ and
$\nu$ satisfies 
\[
\sum_{p} 
\eta_{0}(u_{\lie{n}}^{(p)}) 
(\nu+\rho_{A})(u_{\lie{a}}^{(p)}) 
(\mu_{\delta}+\rhom)(\gamma_{2}(u_{\lie{k}}^{(p)})) 
= \chi_{\Lambda}(z). 
\]
But since 
$\sum_{p} 
\eta_{0}(u_{\lie{n}}^{(p)}) u_{\lie{a}}^{(p)} u_{\lie{k}}^{(p)} 
= \gamma_{1}'(z)$, 
this means that 
\begin{equation}\label{eq:leading term and inf char}
\chi_{(\mu_{\delta}+\rhom, \nu)}(z) = \chi_{\Lambda}(z) 
\quad \Leftrightarrow \quad 
(\mu_{\delta}+\rhom, \nu) \in W \cdot \Lambda. 
\end{equation}
Since $\Lambda$ is regular, $w \Lambda$ ($w \in W$) are all
different. 
Therefore, all the solutions of equation 
\eqref{eq:equation for leading term} are 
\begin{align}
&\phi_{2}^{0}(\delta,\nu,\psi_{1},\psi_{2}; a_{t})(v) 
:= 
t^{\nu+\rho_{A}} \psi_{2} \circ \psi_{1}(v),
\label{eq:leading term of phi_2}
\\
& \psi_{1} \in \Hom_{M}(\tau|_{M},\delta), \
\psi_{2} \in \Hom_{M^{\eta}}(\delta|_{M^{\eta}},\sigma), \
(\mu_{\delta}+\rhom,\nu) \in W \cdot \Lambda. 
\notag
\end{align}
Suppose $\phi_{2}$, $\phi_{2}'$ are two solutions of 
\eqref{eq:diff eq for phi_2}. 
If all the coefficients of 
$\phi_{2}^{0}(\delta,\nu,\psi_{1},\psi_{2}; a_{t})$, 
$\delta \in \widehat{M}$, $\nu \in \lie{a}^{\ast}$, in the
power series expansions of $\phi_{2}$, $\phi_{2}'$ are identical, 
then $\phi_{2}=\phi_{2}'$. 
It follows that the dimension of the solution space of 
\eqref{eq:diff eq for phi_2}, i.e. 
$\dim \stWhLs(\tau)$, is estimated as 
\begin{align}
\dim \stWhLs(\tau) 
\leq 
\sum_{(\delta,\nu) \in \mathcal{A}_{\Lambda}} 
\dim \Hom_{M}(\tau|_{M},\delta)
\dim \Hom_{M^{\eta}}(\delta|_{M^{\eta}},\sigma).
\label{eq:first estimate of K-multiplicity} 
\end{align}

Let 
$\WhL(\tau) 
= \Hom_{K}(V_{\tau}^{K}, \WhL)$. 
By \eqref{eq:first estimate of K-multiplicity}, we have 
\begin{align}
\dim \WhL(\tau) 
&= \sum_{\sigma \in \widehat{M^{\eta}}} 
\dim \stWhLs(\tau) \dim \sigma 
\label{eq:second estimate of K-multiplicity}
\\
& \leq 
\sum_{\sigma \in \widehat{M^{\eta}}} 
\sum_{(\delta,\nu) \in \mathcal{A}_{\Lambda}}
\dim \Hom_{M}(\tau|_{M},\delta)
\dim \Hom_{M^{\eta}}(\delta|_{M^{\eta}},\sigma)
\dim \sigma 
\notag\\
&= 
\sum_{(\delta,\nu) \in \mathcal{A}_{\Lambda}}
\dim \Hom_{M}(\tau|_{M},\delta)
\dim \delta.  
\notag
\end{align}

In the fundamental paper \cite{L}, Lynch gives the dimension of the
space of dual Whittaker vectors of principal series representations. 
His result, together with Theorem~C in \cite{M1}, says that, 
if $(\mu_{\delta}+\rhom, \nu) \in W \cdot \Lambda$, 
then 
\[
\dim \Hom_{\gK}(X_{P}(\delta,\nu), \WhL)
= \# W(G,A) \dim \delta.
\]
Since $\Lambda$ is generic, 
(i) every non-zero element in 
$\Hom_{\gK}(X_{P}(\delta,\nu), \WhL)$ is
injective, and 
(ii) if 
$X_{P}(\delta_{1},\nu_{1})$ and $X_{P}(\delta_{2},\nu_{2})$ are not
equivalent, then for any 
$\Phi_{i} 
\in \Hom_{\gK}(X_{P}(\delta_{i},\nu_{i}), \WhL)$,
$i=1,2$, $\mathrm{Image} \Phi_{1} \cap \mathrm{Image} \Phi_{2} = 0$. 
Then we have 
\begin{align*}
\dim \WhL(\tau) 
& \geq 
\sum_{X_{P}(\delta,\nu) \in \mathcal{H}_{\Lambda}} 
\dim \Hom_{\gK}(X_{P}(\delta,\nu), \WhL)
\dim \Hom_{K}(\tau, X_{P}(\delta,\nu)) 
\\
&= 
\sum_{[(\delta,\nu)] \in W(G,A) \backslash \mathcal{A}_{\Lambda}} 
\# W(G,A) 
\dim \delta 
\times 
\dim \Hom_{M}(\tau|_{M}, \delta) 
\\
& \underset{\eqref{eq:second estimate of K-multiplicity}}{\geq} 
\dim \WhL(\tau). 
\end{align*}
Here, we first used the Frobenius reciprocity 
$\Hom_{K}(\tau, X_{P}(\delta,\nu))
\simeq \Hom_{M}(\tau|_{M}, \delta)$, 
and used the fact that every $W(G,A)$-orbit in
$\mathcal{A}_{\Lambda}$ consists of $\# W(G,A)$ elements, since
$\Lambda$ is regular. 
It follows that every composition factor of 
$\WhL$ 
is a submodule of it. 
In other words, the modules 
$\WhL$, $\stWhLs$ and $\stWhLsmg$ are completely
reducible. 
It also follows that the equality in 
\eqref{eq:first estimate of K-multiplicity} holds.

In order to complete the proof, we recall a result of Wallach's (\cite{W2}). 
Let $dn$ be the Haar measure of $N$ and $w_{0}$ be the longest element
of $W(G,A)$ (with respect to $N$). 
Recall the Jacquet integral 
\begin{align}
J_{\nu} : &C^{\infty}\mbox{-}\Ind_{P}^{G}(\delta \otimes e^{\nu+\rho_{A}}) 
\rightarrow C^{\infty}
(G/M^{\eta}N; (\delta|_{M^{\eta}}) \otimes \eta), 
\label{eq:Jacquet}\\
& 
J_{\nu}(f)(g) 
= 
\int_{N} f(g n w_{0}) \eta(n) dn. 
\notag
\end{align}
Note that $J_{\nu}$ is right $M^{\eta}$-equivariant, since we may
choose $w_{0}$ to commute with $M^{\eta}$. 
Let $\Ind_{P}^{G}(\delta \otimes e^{\nu+\rho_{A}})'$ be the
continuous dual space of 
$C^{\infty}\mbox{-}\Ind_{P}^{G}
(\delta \otimes e^{\nu+\rho_{A}})$, and 
$\Wh_{\eta}^{-\infty}(X_{P}(\delta,\nu))$
be the space of Whittaker vectors in it. 
Let $C^{\infty} \mbox{-}\Ind_{M}^{K}(\delta)$ be the $C^{\infty}$
induced representation of $K$. 
As a $K$-representation, this is isomorphic to 
$C^{\infty}\mbox{-}\Ind_{P}^{G}
(\delta \otimes e^{\nu+\rho_{A}})|_{K}$. 
Let $(\Ind_{M}^{K}(\delta))'$ be the space of all continuous
functionals on $C^{\infty} \mbox{-}\Ind_{M}^{K}(\delta)$, 
which is endowed with the $C^{\infty}$-topology. 
\begin{theorem}[\cite{W2}]\label{theorem:Wallach}
Let $v^{\ast} \in (V_{\delta}^{M})^{\ast}$. 
Then 
$\nu \mapsto \langle v^{\ast}, 
J_{\nu}(\cdot)(e)\rangle_{\delta}$
extends to a weakly holomorphic map of $\lie{a}^{\ast}$ into 
$(\Ind_{M}^{K}(\delta))'$. 
Moreover, for any $\nu \in \lie{a}^{\ast}$, 
\begin{align}
& 
(V_{\delta}^{M})^{\ast}
\ni v^{\ast} 
\mapsto 
(f \mapsto \langle v^{\ast}, J_{\nu}(f)(e)\rangle_{\delta}) 
\in 
\Wh_{\eta}^{-\infty}(X_{P}(\delta,\nu))
\label{eq:continuous Whittaker vector}
\end{align}
is an isomorphism of vector spaces. 
Here, $\langle \ , \ \rangle_{\delta}$ is the pairing of
$V_{\delta}^{M}$ and its dual. 
\end{theorem}
The image of a continuous dual Whittaker vector is characterized
by the moderate growth condition (\cite{W}). 
From this theorem and the map \eqref{eq:Jacquet}, 
we know that there are 
$m_{\delta}(\sigma) 
= \dim \Hom_{M^{\eta}}(\delta|_{M^{\eta}},\sigma)$ copies of 
$X_{P}(\delta,\nu)$ in the socle of $\stWhLsmg$, 
if $X_{P}(\delta, \nu) \in \mathcal{H}_{\Lambda}$. 
As we noted before the theorem, every member of
$\mathcal{H}_{\Lambda}$ is a principal series. 
So the socle of $\stWhLsmg$ is the right hand side
of \eqref{eq:composition series, generic-2}. 
Since $\stWhLsmg$ is completely reducible, the theorem
is shown. 
\end{proof}

\begin{corollary}\label{corollary:K-admissibility}
The standard Whittaker $\brgK$-modules are $K$-admissible 
and they have finite length. 
\end{corollary}
\begin{proof}
If $\Lambda$ is regular, 
the multiplicity of each $K$-type is finite because of 
\eqref{eq:first estimate of K-multiplicity}. 
Such estimate is possible even if $\Lambda$ is not regular. 
The second assertion is clear since these modules admit an
infinitesimal character and $K$-admissible. 
\end{proof}


\section{The group $U(n,1)$ and its irreducible Harish-Chandra
  modules} 
\label{section:U(n,1)}

Up to now very little is known about the properties of 
$\stWhLsmg$ with non-generic $\Lambda$, 
so no smart technique can be used for the analysis of it. 
Therefore, we will choose a group $G$ such that the structure of
Harish-Chandra modules of it is well know and simple (for example
$K$-multiplicity free), 
and we determine the $\brgK$-module structure of 
$\stWhLsmg$ for non-generic $\Lambda$ by direct
calculation. 
Such an example is expected to be a good guide to general cases. 

For such reasons, we assume $G=U(n,1)$ and $\Lambda$ is regular 
integral hereafter.

\subsection{Structure of $U(n,1)$} 
\label{subsection:structure of U(n,1)}

Denote by $E_{ij}$ the standard generators of
$\lie{gl}_{n+1}(\C)$ and define 
$I_{n,1} = \sum_{p=1}^{n} E_{pp} - E_{n+1,n+1}$. 
Let $G = U(n,1)$ be the subgroup of $GL(n+1, \C)$ consisting of
the matrices $g$ satisfying ${}^{t}\bar{g} I_{n,1} g = I_{n,1}$. 
The Lie algebra $\lier{g} = \lie{u}(n,1)$ consists of those matrices
$X \in \lie{gl}_{n+1}(\C)$ which satisfy 
${}^{t}\bar{X} I_{n,1} + I_{n,1} X = O$.
Let $\theta g = I_{n,1} g I_{n,1}$ be a Cartan involution of
$G$. 
The corresponding maximal compact subgroup $K$ of $G$ is 
\[
K = 
\left\{\left.
\begin{pmatrix} 
k & 0 \\ 0 & k_{n+1} 
\end{pmatrix}
\right|
k \in U(n), k_{n+1} \in U(1)
\right\}.
\]
Let $\lier{g}=\lier{k}+\lier{s}$ be the corresponding Cartan
decomposition of $\lier{g}$. 
Then 
\begin{align*}
\{
E_{i,n+1}+E_{n+1,i}, \sqrt{-1}(E_{i,n+1}-E_{n+1,i}) \,|\, 
1 \leq i \leq n \} 
\end{align*}
is a basis of $\lier{s}$. 

Let 
\[
h := E_{n,n+1}+E_{n+1,n} 
\qquad 
\lier{a} 
:= 
\R h, 
\]
and define $f \in \lier{a}^{\ast}$ by $f(h) = 1$. 
Then $\lier{a}$ is a maximal abelian subspace of $\lier{s}$. 
The restricted root system 
$\roots{\lier{g}}{\lier{a}}$ is 
\[
\roots{\lier{g}}{\lier{a}} 
= 
\{\pm f, \pm 2 f\}. 
\]
Choose a positive system 
\[
\proots{\lier{g}}{\lier{a}} 
= 
\{f, 2f\}, 
\]
and denote the corresponding nilpotent subalgebra 
$\sum_{\alpha \in \proots{\lier{g}}{\lier{a}}} (\lier{g})_{\alpha}$ 
by $\lier{n}$. 
One obtains an Iwasawa decomposition 
\begin{align*}
& 
\lier{g} = \lier{k} + \lier{a} + \lier{n}, 
& 
& 
G = K A N, 
\end{align*}
where $A = \exp \lier{a}$ and $N = \exp \lier{n}$. 
Let 
\begin{equation}\label{eq:basis of n}
\begin{split}
& 
X_{i} 
:= 
E_{in}-E_{ni}-E_{i,n+1}-E_{n+1,i}
\qquad 
(1 \leq i \leq n-1), 
\\
& 
Y_{i} 
:= 
\I(E_{in}+E_{ni}-E_{i,n+1}+E_{n+1,i})
\qquad 
(1 \leq i \leq n-1), 
\\
&
Z
:=\I(E_{nn}-E_{n+1,n+1}-E_{n,n+1}+E_{n+1,n}). 
\end{split}
\end{equation}
Then 
$\{X_{i}, Y_{i} \,|\, 1 \leq i \leq n-1\}$ is a basis of
$(\lier{g})_{f}$, 
and 
$\{Z\}$ is a basis of 
$(\lier{g})_{2f}$. 

In our $U(n,1)$ case, 
$M$ is isomorphic to $U(n-1) \times U(1)$. 
It acts on the space of non-degenerate unitary characters of $N$ by 
$\eta \mapsto \eta^{m}(n) := \eta(m^{-1}nm)$, $m \in M$. 
Therefore, we may choose a manageable unitary character when we
calculate Whittaker modules. 
We use the non-degenerate character $\eta$ defined by 
\begin{equation}\label{eq:definition of psi}
\begin{array}{ll}
\eta(X_{i}) = 0, \quad i=1,\dots,n-1,
& \eta(Y_{i}) = 0, \quad i=1,\dots,n-2,
\\
\eta(Y_{n-1}) = \I \xi, \quad \xi > 0, 
& \eta(Z) = 0. 
\end{array}
\end{equation}
It is easy to see that $M^{\eta}$ is isomorphic to 
$U(n-2) \times U(1)$.

\subsection{Classification of irreducible Harish-Chandra modules}
\label{subsection:classification of (g,K)-modules}

We review the 
classification of irreducible Harish-Chandra modules of $G=U(n,1)$
with regular integral infinitesimal character. 
For details, see \cite{C}, \cite{Kr} for example. 
We use the notation $\pi_{i,j}$, $\overline{\pi}_{i,j}$ etc in
\cite{C}. 

There are two conjugacy classes of Cartan subgroups in $G$,
one is compact and the other is maximally split.  
Let $H_{c}$ be the compact Cartan subgroup consisting of diagonal
matrices and let $\lie{h}_{c}$ be its Lie algebra. 
Define a basis  $\{\epsilon_{i}\,|\, i=1,\dots,n+1\}$ of $\lie{h}_{c}^{\ast}$ 
by $\epsilon_{i}(E_{jj}) = \delta_{ij}$. 
Choose a maximally split Cartan subgroup $H_{s}:=(H_{c}\cap M) A$. 
The complexified Lie algebra $\lie{h}_{s}$ of it is the linear span of 
$E_{ii}$ ($i=1,\dots,n-1$), $E_{nn}+E_{n+1,n+1}$ and $h$. 
Define $\tilde{\epsilon}_{n} \in (\lie{h}_{c} \cap \lie{m})^{\ast}$ by 
$\tilde{\epsilon}_{n}(E_{jj}) = 0$ for $j=1,\dots,n-1$ and
$\tilde{\epsilon}_{n}(E_{nn} + E_{n+1,n+1}) = 1$. 
Then $\epsilon_{i}$, $i=1,\dots,n-1$, $\tilde{\epsilon}_{n}$ and $f$ is a basis of
$\lie{h}_{s}^{\ast}$.

Consider the irreducible Harish-Chandra modules with the
regular integral infinitesimal character $\Lambda$, 
which is conjugate to 
\begin{equation}\label{eq:Lambda}
\sum_{p=1}^{n+1} \Lambda_{p} \epsilon_{p} \in \lie{h}_{c}^{\ast}, 
\quad 
\Lambda_{p} \in \Z+n/2, 
\quad 
\Lambda_{1} > \Lambda_{2} > \dots > \Lambda_{n+1}. 
\end{equation}
There are $n+1$ inequivalent discrete series representations
$\pi_{i}$, $i=0,\dots,n$, whose Harish-Chandra parameters are 
\begin{align*}
\sum_{p=1}^{i} \Lambda_{p} \epsilon_{p} 
+ \sum_{p=i+1}^{n} \Lambda_{p+1} \epsilon_{p} 
+ \Lambda_{i+1} \epsilon_{n+1}, 
\end{align*}
respectively. 
$\pi_{i}$ is also denoted by $\overline{\pi}_{i,n+1-i}$. 

For $i=0,\dots,n-1$ and $j=1,\dots,n-i$, 
define $\mu_{i,j} \in (\lie{h}_{c} \cap \lie{m})^{\ast}$ and 
$\nu_{i,j} \in \lie{a}^{\ast}$ by 
\begin{equation}\label{eq:M-pseudocharacter}
\begin{split}
\mu_{i,j} := 
& \sum_{p=1}^{i} \Lambda_{p} \epsilon_{p} 
+ \sum_{p=i+1}^{n-j} \Lambda_{p+1} \epsilon_{p} 
+ \sum_{p=n-j+1}^{n-1} \Lambda_{p+2} \epsilon_{p} 
- \rhom 
\\
& \qquad +(\Lambda_{i+1}+\Lambda_{n-j+2}) \tilde{\epsilon}_{n} 
\\
\nu_{i,j} :=& 
(\Lambda_{i+1}-\Lambda_{n-j+2})f, 
\end{split}
\end{equation}
where 
$\rhom :=\frac{1}{2} \sum_{p=1}^{n-1} (n-2p) \epsilon_{p}$. 
Let $\delta_{i,j}$ be the irreducible representation of $M$ with
the highest weight $\mu_{i,j}$, 
and let $\pi_{i,j}:=X_{P}(\delta_{i,j},\nu_{i,j})$. 
Then $\pi_{i,j}$ has the unique irreducible quotient, which we denote
by $\overline{\pi}_{i,j}$. 
\begin{theorem}
The irreducible Harish-Chandra modules of $U(n,1)$ with the regular
integral infinitesimal character $\Lambda$ are parametrized, 
up to $K$-conjugacy, 
by the set 
$\{\overline{\pi}_{i,j}\,|\, i=0,\dots,n, j=1,\dots,n+1-i\}$. 
\end{theorem}

The $K$-type structure of $\overline{\pi}_{i,j}$ is explicitly known. 
To state the theorem, let $\Lambda_{0} := \infty$ and 
$\Lambda_{n+2}:= -\infty$. 
As is explained in \S~\ref{section:introduction}, 
write $\vert v \vert := \sum_{i=1}^{n+1} v_{i}$ 
for an element 
$v = \sum_{i=1}^{n+1} v_{i} \epsilon_{i} \in \lie{h}_{c}^{\ast}$. 
\begin{theorem}[\cite{Kr}]
\label{theorem:K-spectra} 
For $i=0,\dots,n$ and $j=1,\dots,n+1-i$, 
the $K$-spectrum $\widehat{K}(\overline{\pi}_{i,j})$ is 
\begin{equation}\label{eq:K-spectrum of pi_ij}
\begin{split}
\{(\tau_{\lambda}, V_{\lambda}^{K}) \,|\, 
\, & 
\Lambda_{p-1} -n/2+p-1 
\geq \lambda_{p} \geq
\Lambda_{p}-n/2+p, 
\enskip
p=1,\dots,i; 
\\
& 
\Lambda_{p} -n/2+p-1 
\geq \lambda_{p} \geq
\Lambda_{p+1}-n/2+p, 
\enskip 
p=i+1,\dots,n-j+1; 
\\
& 
\Lambda_{p+1} -n/2+p-1 
\geq \lambda_{p} \geq
\Lambda_{p+2}-n/2+p, 
\enskip 
p=n-j+2,\dots,n;
\\
& \vert \lambda \vert = \vert \Lambda \vert
\}, 
\end{split}
\end{equation}
and each $K$-type occurs in $\overline{\pi}_{i,j}$ with multiplicity one. 
\end{theorem}

In order to state the composition series, we use diagrammatic
expression. 
\begin{definition}\label{definition:diagram of composition series} 
Suppose $A_{1}, A_{2}$ are distinct composition factors of a 
$\brgK$-module $V$. 
If there exist elements $\{v_{i}\} \subset A_{1}$ and 
$\{X_{i}\} \subset \lie{g}$ such that $\sum_{i} X_{i} v_{i}$ is
non-zero and contained in $A_{2}$, 
then we connect $A_{1}$ and $A_{2}$ by the arrow 
$A_{1} \rightarrow A_{2}$. 
\end{definition}

\begin{theorem}[\cite{Kr}, \cite{C}]\label{composition series of PS} 
The composition series of $\pi_{i,j}$, 
$i=0,\dots,n-1$, $j=1,\dots,n-i$, is 
\[
\pi_{i,j} 
= 
X_{P}(\delta_{i,j},\nu_{i,j}) 
\quad \simeq \quad 
\begin{xy}
(0,7)*{\overline{\pi}_{i,j}}="A_{11}",
(-7,0)*{\overline{\pi}_{i,j+1}}="A_{21}",
(7,0)*{\overline{\pi}_{i+1,j}}="A_{22}",
(0,-7)*{\overline{\pi}_{i+1,j+1}}="A_{31}",
\ar "A_{11}";"A_{21}"
\ar "A_{11}";"A_{22}"
\ar "A_{21}";"A_{31}"
\ar "A_{22}";"A_{31}"
\end{xy}
\]
If $i+j=n$, the factor $\overline{\pi}_{i+1,j+1}$ does not appear. 
\end{theorem}


\section{Composition factors of $\stWhLsmg$}
\label{section:composition factors}

In this section we first determine the submodules of
$\stWhLsmg$. 
For this purpose, we need some results on the Whittaker models.

\subsection{Whittaker models}
\label{subsection:Whittaker models}
Let $(\pi, V)$ be an irreducible Harish-Chandra module. 
A realization of $(\pi, V)$ as a submodule of 
$C^{\infty}(G/N;\eta)_{K}$ 
is called a Whittaker model of $(\pi, V)$. 
For a Harish-Chandra module $V$, 
let $V_{\infty}$ be its $C^{\infty}$-globalization. 
As in \S\ref{section:generic case}, let $\Wh_{\eta}^{-\infty}(V)$ 
be the space of Whittaker vectors in the continuous dual space of
$V_{\infty}$. 
Note that the image of an element of $\Wh_{\eta}^{-\infty}(V)$ is
characterized by the moderate growth condition. 
The next theorem tells us which irreducible $\brgK$-module
can be a submodule of $C^{\infty}(G/N;\eta)_{K}$. 
\begin{theorem}[\cite{M1}, \cite{M2}]\label{theorem:matu-1}
Let $V$ be a Harish-Chandra module. 
\begin{enumerate}
\item
$V$ has a non-trivial Whittaker model 
if and only if the Gelfand-Kirillov dimension $\Dim V$ of $V$ is equal
to $\dim N$. 
\item (Casselman) 
$V \to \Wh_{\eta}^{-\infty}(V)$ is an exact functor. 
\end{enumerate}
\end{theorem}

The Gelfand-Kirillov dimensions of the irreducible 
modules $\overline{\pi}_{i,j}$ are 
\begin{align*}
& \Dim \overline{\pi}_{0,1} = 0, 
\\
& \Dim \overline{\pi}_{i,1} = \Dim  \overline{\pi}_{0,j} = n, 
& 
& i=1,\dots,n, \quad j=2,\dots, n+1, 
\\
& \Dim \overline{\pi}_{i,j} = 2n-1 = \dim N, 
&
& i=1,\dots,n-1, \quad j=2,\dots,n+1-i. 
\end{align*}
See \cite{C} for example. 
Therefore, an irreducible submodule of $\stWhLsmg$ is isomorphic to 
one of $\overline{\pi}_{i,j}$, 
$i=1,\dots,n$, $j=2,\dots,n+1-i$. 

\subsection{Unique simple submodule} 
\label{subsection:Unique simple submodule}
Let $(\pi, V)$ be an irreducible Harish-Chandra module with 
$\Dim V = \dim N$. 
Suppose that it is an composition factor of some principal series
representation $X_{P}(\delta,\nu)$. 
By Theorem~\ref{theorem:matu-1}(2) and Theorem~\ref{theorem:Wallach}, 
every continuous embedding of $V_{\infty}$ into 
$C^{\infty}(G/N;\eta)$ is a composition of 
(i) a realization of $V_{\infty}$ as a subquotient of 
$C^{\infty}$-$\Ind_{P}^{G}(\delta \otimes e^{\nu+\rho_{A}})$  
and (ii) a Jacquet integral. 
Since a Jacquet integral is right $M^{\eta}$-equivariant, 
$V$ can be a submodule of $\stWhLsmg$ only if 
$\sigma \subset \delta|_{M^{\eta}}$. 
Let $\{X_{P}(\delta_{p},\nu_{p})\, | \, p=1,\dots,k\}$ 
be the set of principal series representations which contain 
$(\pi,V)$ as a subquotient. 
If $(\pi,V)$ is a submodule of $\stWhLsmg$, then by the
above discussion 
$\sigma 
\subset \cap_{p=1}^{k} \delta_{p}|_{M^{\eta}}$, 
i.e. $\sigma$ is a submodule of $\delta_{i}|_{M^{\eta}}$ for every 
$p=1,\dots,k$. 

Conversely, for $\sigma \in \widehat{M^{\eta}}$, 
suppose that there exists a principal series 
$X_{P}(\delta,\nu) \in \mathcal{H}_{\Lambda}$ such that 
$\sigma \subset \delta|_{M^{\eta}}$. 
Then by Theorem~\ref{theorem:Wallach}, the intersection of the image
of the Jacquet integral \eqref{eq:Jacquet} and $\stWhLsmg$ 
is non-zero. 
Especially, $\stWhLsmg$ is non-zero. 

\begin{proposition}
\label{proposition:condition for sigma, integral} 
Suppose the regular infinitesimal character $\Lambda$ is
integral. 
The irreducible module  
$\overline{\pi}_{i,j}$, $i=1,\dots,n-1$, $j=2,\dots,n+1-i$, is a
submodule of $\stWhLsmg$ 
if and only if the highest weight 
$\gamma = (\gamma_{1},\dots,\gamma_{n-2}; \gamma_{n-1})$ of the
irreducible representation $\sigma$ of 
$M^{\eta} \simeq U(n-2) \times U(1)$
satisfies 
\begin{equation}\label{eq:condition for sigma} 
\begin{cases}
&\Lambda_{p}-n/2+p 
\geq \gamma_{p} \geq 
\Lambda_{p+1}-n/2+p+1, 
\quad 
p=1,\dots,i-1, 
\\
&\Lambda_{p+1}-n/2+p 
\geq \gamma_{p} \geq 
\Lambda_{p+2}-n/2+p+1, 
\quad 
p=i,\dots,n-j, 
\\
&\Lambda_{p+2}-n/2+p 
\geq \gamma_{p} \geq 
\Lambda_{p+3}-n/2+p+1, 
\quad 
p=n-j+1,\dots,n-2, 
\\
& \gamma_{n-1} = \vert \Lambda \vert 
- \sum_{p=1}^{n-2} \gamma_{p}.  
\end{cases}
\end{equation}
Especially, $\stWhLsmg$ is non-zero 
if and only if the highest weight of
$\sigma$ satisfies the condition \eqref{eq:condition for sigma} for
some $i, j$. 
In this case, $\overline{\pi}_{i,j}$ is the unique simple submodule of
$\stWhLsmg$. 
\end{proposition}
\begin{proof}
By Theorem~\ref{composition series of PS}, 
$\overline{\pi}_{i,j}$, $i=1,\dots,n$, $j=2,\dots,n+1-i$, 
is a composition factor of the principal series
$\pi_{k,l}$ if and only if 
$(k,l) = (i,j)$ (only when $i+j \leq n$), $(i,j-1), (i-1,j)$ or
$(i-1,j-1)$. 
Therefore, if $\overline{\pi}_{i,j}$ is a submodule of
$\stWhLsmg$, then 
$\sigma 
\subset 
\delta_{i,j}|_{M^{\eta}} \cap \delta_{i,j-1}|_{M^{\eta}} 
\cap \delta_{i-1,j}|_{M^{\eta}} \cap \delta_{i-1,j-1}|_{M^{\eta}}$. 
Conversely, if $\sigma$ satisfies this condition, 
then $\stWhLsmg$ is non-zero, as is stated before this proposition. 

Recall the branching rule for the restriction $U(m)$ to $U(m-1)$. 
For an irreducible representation $\delta_{\mu}$ of $U(m)$ with the
highest weight $\mu=(\mu_{1},\dots,\mu_{m})$, 
the restriction $\delta_{\mu}|_{U(m-1)}$ is a direct sum of
$\sigma_{\gamma} \in U(m-1)$, with 
\begin{align*}
& \gamma = (\gamma_{1},\dots,\gamma_{m-1}), 
& 
& \mu_{p} \geq \gamma_{p} \geq \mu_{p+1}, 
\quad p=1,\dots,m-1, 
\quad \gamma_{p} \in \Z. 
\end{align*}
It follows that the restriction 
$\delta_{k,l} \in \widehat{M} \simeq \widehat{U(n-1)} \times
\widehat{U(1)}$ to $M^{\eta} \simeq U(n-2) \times U(1)$ is a direct
sum of $\sigma_{\gamma} \in \widehat{U(n-2)} \times \widehat{U(1)}$,
whose highest weight 
$\gamma = (\gamma_{1},\dots,\gamma_{n-2};\gamma_{n-1})$ satisfies 
\[
\begin{cases}
\Lambda_{p}-n/2+p 
\geq \gamma_{p} \geq 
\Lambda_{p+1}-n/2+p+1, 
\quad 
p=1,\dots,k-1, 
\\
\Lambda_{k}-n/2+k 
\geq \gamma_{k} \geq 
\Lambda_{k+2}-n/2+k+1, 
\\
\Lambda_{p+1}-n/2+p 
\geq \gamma_{p} \geq 
\Lambda_{p+2}-n/2+p+1, 
\quad 
p=k+1,\dots,n-l-1, 
\\
\Lambda_{n-l+1}+n/2-l 
\geq \gamma_{n-l} \geq 
\Lambda_{n-l+3}+n/2-l+1, 
\\
\Lambda_{p+2}-n/2+p 
\geq \gamma_{p} \geq 
\Lambda_{p+3}-n/2+p+1, 
\quad 
p=n-l+1,\dots,n-2, 
\\
\gamma_{n-1} = \vert \Lambda \vert 
- \sum_{p=1}^{n-2} \gamma_{p},  
\end{cases}
\]
if $k+l \leq n-1$. 
The last condition for $\gamma_{n-1}$ is obtained from the action
of the center of $G$. 
Therefore, if $i+j \leq n$, then $\sigma \in \widehat{M^{\eta}}$
satisfies 
$\sigma 
\subset 
\delta_{i,j}|_{M^{\eta}} \cap \delta_{i,j-1}|_{M^{\eta}} 
\cap \delta_{i-1,j}|_{M^{\eta}} \cap \delta_{i-1,j-1}|_{M^{\eta}}$
if and only if the highest weight $\gamma$ of $\sigma$ satisfies 
\eqref{eq:condition for sigma}. 
This proves the ``only if'' part of proposition for the case 
$i+j\leq n$. 
The case $i+j=n+1$ is shown analogously. 

We will show that the condition is sufficient and that the
multiplicity in the socle is one. 

Let $(\widehat{M^{\eta}})_{i,j}$ be the set of 
$\sigma \in \widehat{M^{\eta}}$ whose highest weight $\gamma$
satisfies the condition \eqref{eq:condition for sigma}. 
Then it is easy to see that 
$(\widehat{M^{\eta}})_{i,j} 
\cap 
(\widehat{M^{\eta}})_{k,l} 
= \emptyset$ if $(i,j) \not= (k,l)$. 
It follows that if 
$\sigma \in (\widehat{M^{\eta}})_{i,j}$, 
then every irreducible factor in the socle of 
$\stWhLsmg$ is isomorphic to $\overline{\pi}_{i,j}$. 
Let $m_{\sigma}$ be the multiplicity of such factors. 
Then 
$\dim \Wh_{\eta}^{-\infty}(\overline{\pi}_{i,j}) 
= \sum_{\sigma \in (\widehat{M^{\eta}})_{i,j}} 
m_{\sigma} \dim \sigma$. 
By Theorems~\ref{composition series of PS} and 
\ref{theorem:matu-1}, we have 
\[
\dim \Wh_{\eta}^{-\infty}(\pi_{i,j})
= 
\sum_{a,b=0,1}
\dim \Wh_{\eta}^{-\infty}(\overline{\pi}_{i+a,j+b}) 
=
\sum_{a,b=0,1} 
\sum_{\sigma \in (\widehat{M^{\eta}})_{i+a,j+b}} 
m_{\sigma} \dim \sigma. 
\]
On the other hand, it is easy to see from
Theorem~\ref{theorem:Wallach} and \eqref{eq:condition for sigma} that 
\[
\dim \Wh_{\eta}^{-\infty}(\pi_{i,j})
= \dim \delta_{i,j} 
= 
\sum_{a,b=0,1} 
\sum_{\sigma \in (\widehat{M^{\eta}})_{i+a,j+b}} 
\dim \sigma. 
\]
Since $m_{\sigma} \geq 1$ for any 
$\sigma \in \cup_{a,b=0,1} (\widehat{M^{\eta}})_{i+a,j+b}$, 
they are all one. 
This completes the proof of proposition. 
\end{proof}

\subsection{Composition factors}\label{subsection:composition factors} 

Hereafter, we denote $\stWhLsmg$ by 
$\stWhLgmg$ if the highest weight of $\sigma$ is
  $\gamma$. 
We also denote by $\sigma_{\gamma}$ the irreducible representation of
$M^{\eta}$ whose highest weight is $\gamma$. 
We first determine the irreducible representations appearing in the
composition series of $\stWhLgmg$. 
\begin{proposition}\label{proposition:composition factors-1} 
Suppose that $\Lambda$ is regular integral and 
that $\gamma$ satisfies \eqref{eq:condition for sigma}, 
so $\overline{\pi}_{i,j}$ is the unique simple submodule of
$\stWhLgmg$. 
In this case, an irreducible module $\overline{\pi}_{k,l}$ is a
composition factor of $\stWhLgmg$ only if 
$(k,l)=(i+a,j+b)$ with $a=0, \pm 1$ and $b=0, \pm 1$.  
\end{proposition}
\begin{proof}
We have seen in \S~\ref{section:generic case} that, 
if $\overline{\pi}_{k,l}$ is a composition factor of
$\stWhLgmg$, each $K$-type of it must contain the
representation $\sigma_{\gamma}$. 
By \eqref{eq:K-spectrum of pi_ij} and \eqref{eq:condition for sigma}, 
this is possible if and only if 
$(k,l)=(i+a,j+b)$ with $a=0, \pm 1$ and $b=0, \pm 1$.  
\end{proof}

\begin{proposition}\label{proposition:at least one}
Suppose that $\gamma$ satisfies \eqref{eq:condition for sigma} 
and $(i',j') = (i+a,j+b)$, $a,b=0,\pm 1$. 
Then the multiplicity of $\overline{\pi}_{i',j'}$ in
$\stWhLgmg$ is at least one. 
\end{proposition}
\begin{proof}
By Theorem~\ref{composition series of PS}, $\overline{\pi}_{i',j'}$
is a composition factor of 
$\pi_{k,l}=X_{P}(\delta_{k,l},\nu_{k,l})$, 
with $(k,l)=(i+a,j+b)$,  $a,b=0,-1$. 
Consider the principal series representation 
$X_{P}(\delta_{k,l},\nu)$ with $\nu \in \lie{a}^{\ast}$. 
If $\nu$ is generic, then
$\stWhmg{(\mu_{k,l}+\rhom,\nu)}{\gamma}$ is isomorphic to 
$X_{P}(\delta_{k,l},\nu)$ by Theorem~\ref{theorem:generic case}. 
Here we used the fact $m_{\delta_{k,l}}(\sigma_{\gamma}) = 1$. 
Choose a $K$-type $\tau$ of $\overline{\pi}_{i',j'}$. 
This is also a $K$-type of $X_{P}(\delta_{k,l},\nu)$. 
By Frobenius reciprocity, the multiplicity of $\tau$ in
$X_{P}(\delta_{k,l},\nu)$ is one. 
Therefore, the space of moderately growing solutions of 
\eqref{eq:diff eq for phi_2}, with $\Lambda$ replaced by 
$(\mu_{k,l}+\rhom,\nu)$ and $\sigma$ by $\sigma_{\gamma}$, 
is one dimensional. 
Let $f_{\nu}$ be a non-zero moderately growing solution. 
Then by Theorem~\ref{theorem:Wallach}, this function is expressed by
the Jacquet integral and it is holomorphic in $\nu$. 
Suppose the order of zero of $f_{\nu}$ at $\nu=\nu_{k,l}$ is $m$. 
Let $g_{\nu} :=f_{\nu}/(\nu-\nu_{k,l})^{m}$. 
Then $g_{\nu_{k,l}}$ is non-zero. 
It satisfies the equation 
\eqref{eq:diff eq for phi_2} (with $\sigma$ replaced by
$\sigma_{\gamma}$) and grows moderately at the infinity, 
so it is an element of $\stWhLgmg$. 

We have proved that, for every $K$-type $\tau$ of
$\overline{\pi}_{i',j'}$, the multiplicity of $\tau$ in
$\stWhLgmg$ is at least one. 
Since 
$\widehat{K}(\overline{\pi}_{a,b}) 
\cap
\widehat{K}(\overline{\pi}_{a',b'}) 
= \emptyset$ if $(a,b)\not=(a',b')$, 
the multiplicity of $\overline{\pi}_{i',j'}$ in $\stWhLgmg$ is at
least one. 
\end{proof}


\section{Determination of the composition series}
\label{section:determination of the composition series}

In this section, we determine the composition series of
$\stWhLgmg$ in the case when $\Lambda$ is integral. 
For this purpose, we need to write down the actions of $Z(\lie{g})$
and $\lie{s}$ on this space explicitly. 
The former is achieved by the determinant type central element of
$U(\lie{gl}_{n+1})$, and the latter by the $K$-type shift operators.

\subsection{Shift operators}\label{subsection:shift operators} 
We review the $K$-type shift operators briefly. 
Choose a $K$-type $(\tau_{\lambda},V_{\lambda}^{K})$ of
$\stWhLsmg$, whose highest weight is $\lambda$. 
Let $\psi$ be an element of 
$\Hom_{K}(V_{\lambda}^{K}, \stWhLsmg)$. 
For $v \in V_{\lambda}^{K}$ and $X \in \lie{s}$, define 
\[
\tilde{\psi}(v \otimes X)(g) 
:= L_{X} (\psi(v))(g). 
\]
Then it is easy to see that $\tilde{\psi}$ is an element of 
$\Hom_{K}(V_{\lambda}^{K} \otimes \lie{s}, \stWhLsmg)$. 
Here, we regard $\lie{s}$ as a representation of $K$ by the adjoint
action $\Ad$. 
Denote by $\Delta_{\lie{s}}$ the set of weights on $\lie{s}$ with
respect to a fixed Cartan subalgebra of $\lie{k}$. 
In our case, the irreducible decomposition of 
$V_{\lambda}^{K} \otimes \lie{s}$ is 
$\oplus_{\alpha \in \Delta_{\lie{s}}}
m(\alpha)\, V_{\lambda+\alpha}^{K}$, 
$m(\alpha) = 0$ or $1$. 
When $m(\alpha) = 1$, let $\iota_{\alpha}$ be the embedding of
$V_{\lambda+\alpha}^{K}$ into 
$V_{\lambda}^{K} \otimes \lie{s}$. 
Define 
\[
\tilde{\psi}_{\alpha}(v_{\alpha})(g) 
:= 
\tilde{\psi}(\iota_{\alpha}(v_{\alpha}))(g), 
\quad v_{\alpha} \in V_{\lambda+\alpha}^{K}. 
\]
Then $\tilde{\psi}_{\alpha}$ is an element of 
$\Hom_{K}(V_{\lambda+\alpha}^{K}, \stWhLsmg)$, 
and the correspondence $\psi \mapsto \tilde{\psi}_{\alpha}$ is a
$K$-type shift in $\stWhLsmg$ coming from the
$\lie{s}$-action.

\subsection{Gelfand-Tsetlin basis} 
\label{subsection:Gelfand-Tsetlin basis} 
In order to write down the $K$-type shift operators explicitly, 
we realize the space 
$\Hom_{M^{\eta}}(V_{\tau}^{K},V_{\sigma}^{M^{\eta}})$ 
by using the Gelfand-Tsetlin basis (\cite{GT}). 

\begin{definition}\label{definition:GT}
Let $\lambda = (\lambda_{1}, \dots, \lambda_{n})$ be a dominant
integral weight of $U(n)$. 
A {\it ($\lambda$-)Gelfand-Tsetlin pattern} is a set of vectors 
$Q=(\vect{q}_{1}, \dots, \vect{q}_{n})$ such that 
\begin{enumerate}
\item
$\vect{q}_i=(q_{1,i}, q_{2,i}, \dots, q_{i,i})$. 
\item
The numbers $q_{i,j}$ are all integers. 
\item
$q_{i,j+1} \geq q_{i,j} \geq q_{i+1,j+1}$, 
for any $i=1, \dots, j$.  
\item
$q_{i,n} = \lambda_{n}$, $i=1,\dots,n$. 
\end{enumerate}
The set of all $\lambda$-Gelfand-Tsetlin patterns is denoted by
$GT(\lambda)$. 
\end{definition}
\begin{theorem}[\cite{GT}] 
For a dominant integral weight $\lambda$ of $U(n)$, 
let $(\tau_{\lambda}, V_{\lambda}^{U(n)})$ be the irreducible
representation of $U(n)$ with the highest weight $\lambda$. 
Then $GT(\lambda)$ is identified with a basis of 
$(\tau_{\lambda}, V_{\lambda}^{U(n)})$. 
\end{theorem}
The action of elements $E_{ij} \in \lie{gl}(n,\C)$ is expressed as
follows. 
Let $l_{i,j} := q_{i,j} - i$ 
and $\vert \vect{q}_{j} \vert := \sum_{i=1}^{j} q_{i,j}$. 
Let $\sigma_{i,j}^{\pm}$ be the shift operators on $GT(\lambda)$, 
sending $\vect{q}_{j}$ to 
$\vect{q}_{j} + (0, \dots, \overset{i}{\pm 1}, 0 ,\dots, 0)$. 
Define $a_{i,j}(Q)$ and $b_{i,j}(Q)$ by 
\begin{align} 
a_{i,j}(Q) 
&= 
\left\vert 
\frac{\prod_{k=1}^{j+1} (l_{k,j+1} - l_{i,j}) 
\prod_{k=1}^{j-1} (l_{k,j-1} - l_{i,j} - 1)}
{\prod_{\genfrac{}{}{0pt}{}{k=1}{k \not= i}}^{j} 
(l_{k,j}-l_{i,j}) (l_{k,j}-l_{i,j}-1)}
\right\vert^{1/2}, 
& 
b_{i,j}(Q) 
&= a_{i,j}(\sigma_{i,j}^{-} Q). 
\label{eq:coefficients of GT}
\end{align}

\begin{theorem}[\cite{GT}]
\label{theorem:GT}
For $Q \in GT(\lambda)$, the action of the Lie algebra is given by 
\begin{align*}
&\tau_\lambda(E_{j,j+1}) Q
=
\sum_{i=1}^{j} a_{i,j}(Q) \sigma_{i,j}^{+} Q, 
&
&\tau_\lambda(E_{j+1,j}) Q
=
\sum_{i=1}^{j} b_{i,j}(Q) \sigma_{i,j}^{-} Q, 
\\
&\tau_{\lambda}(E_{jj}) Q 
= (\vert \vect{q}_{j} \vert - \vert \vect{q}_{j-1} \vert) Q. 
\end{align*}
\end{theorem}

\begin{remark}\label{remark:GT restriction}
The Gelfand-Tsetlin basis is compatible with the restriction to
smaller unitary groups $U(k)$, $k=1,\dots,n-1$. 
More precisely, the restriction of $\tau_{\lambda}$ to $U(n-1)$ is
multiplicity free, and the highest weights of the irreducible
representation appearing in $\tau_{\lambda}|_{U(n-1)}$ are the above 
$\vect{q}_{n-1}$'s. 
\end{remark}

\begin{remark}\label{remark:GT of dual}
The highest weight $\lambda^{\ast}$ of the contragredient representation 
$(\tau_{\lambda^{\ast}}, V_{\lambda^{\ast}}^{K})$ of 
$(\tau_{\lambda}, V_{\lambda}^{K})$ is 
$\lambda^{\ast} = (-\lambda_{n}, \dots, -\lambda_{1})$. 
In this case, 
$Q^{\ast} := (\vect{q}_{1}^{\ast}, \dots, \vect{q}_{n}^{\ast}) 
\in GT(\lambda^{\ast})$, 
$\vect{q}_{i}^{\ast} := (-q_{i,i}, \dots, -q_{1,i})$ is dual to
$Q \in GT(\lambda)$. 
\end{remark}

\subsection{Explicit formulas of shift operators}
\label{subsection:explicit shift operators}
In \S\ref{section:generic case}, we identified 
an element 
$\phi_{1} \in \Hom_{K}(V_{\lambda}^{K},\stWhLg)$ 
with a function 
$\phi_{2} \in \stWhLg(\tau_{\lambda})$ 
(for the definition of this space, see \eqref{eq:tau isotypic in I}). 
The space 
$\Hom_{M^{\eta}}
(V_{\lambda}^{K},V_{\gamma}^{M^{\eta}})$ is
isomorphic to 
$(V_{\lambda^{\ast}}^{K} \otimes 
V_{\gamma}^{M^{\eta}})^{M^{\eta}}$, 
the space of $M^{\eta}$-invariants in 
$V_{\lambda^{\ast}}^{K} \otimes V_{\gamma}^{M^{\eta}}$. 
By Remark~\ref{remark:GT restriction}, a basis of this space
is identified with the ``partial Gelfand-Tsetlin patterns'' 
\begin{align}
GT((\lambda/\gamma)^{\ast}) 
:= 
\{Q = &(\vect{q}_{n-2}, \vect{q}_{n-1}, \vect{q}_{n}) \,|\, 
\label{eq:partial GT basis}
\\
&\vect{q}_{n-2} = \gamma^{\ast}, \vect{q}_{n} = \lambda^{\ast}; \, 
\mbox{ satisfies Definition~\ref{definition:GT} (1)--(3)}\}. 
\notag
\end{align}
The correspondence is given by 
\begin{align}
& GT((\lambda/\gamma)^{\ast}) 
\ni Q \mapsto 
\langle\langle \ast, Q\rangle\rangle_{\lambda} 
\in 
\Hom_{M^{\eta}}(V_{\lambda}^{K},V_{\gamma}^{M^{\eta}}), 
\mbox{ where}
\notag\\
& \langle\langle Q',Q \rangle\rangle_{\lambda} 
= 
\begin{cases} 
0 \quad 
\mbox{if} \ 
\vect{q}_{n-1}'\not=\vect{q}_{n-1}^{\ast} 
\ \mbox{or} \ 
\vect{q}_{n-2}'\not=\gamma, 
\\
(\vect{q}_{1}',\dots,\vect{q}_{n-2}') \in GT(\gamma)  
\quad \mbox{if} \ 
\vect{q}_{n-1}'=\vect{q}_{n-1}^{\ast} 
\ \mbox{and} \ 
\vect{q}_{n-2}'=\gamma, 
\end{cases} 
\label{eq:partial GT <-> Hom}
\\
& \mbox{for} \ 
Q' = (\vect{q}_{1}',\dots,\vect{q}_{n}') \in GT(\lambda). 
\notag
\end{align}
The action of $\lie{k}$ on $GT((\lambda/\gamma)^{\ast})$ is given by
$\langle\langle Q', 
\tau_{\lambda^{\ast}}(\cdot) Q 
\rangle\rangle_{\lambda}
= - 
\langle\langle \tau_{\lambda}(\cdot) Q', Q 
\rangle\rangle_{\lambda}$. 

Let $V_{(\lambda/\gamma)^{\ast}}^{K/M^{\eta}}$ be the vector space
spanned by $GT((\lambda/\gamma)^{\ast})$. 
Then we can identify 
$\phi(a) 
\in 
C^{\infty}(A \rightarrow V_{(\lambda/\gamma)^{\ast}}^{K/M^{\eta}})$ 
with $\phi_{2} \in \stWhLg(\tau_{\lambda})$ 
via 
\[
\phi_{2}(a)(v) 
= 
\langle\langle v, \phi(a) \rangle\rangle_{\lambda}, 
\quad 
v \in V_{\lambda}^{K}, a \in A. 
\]
We introduce a coordinate system on $A$ defined by 
\begin{equation}\label{eq:coordinate on A} 
\R_{>0} \ni t \mapsto 
a_{t}:= \exp(\log(\xi/t) h) \in A. 
\end{equation}
Then the action of $\lie{g}$ on $\phi$ is given by 
\begin{equation}\label{eq:action of Lie algebra elements}
\begin{split}
& h \cdot \phi(a_{t}) = \theta \phi(a_{t}), \quad 
\theta := t \frac{d}{dt}, 
\qquad 
Y_{n-1} \cdot \phi(a_{t}) 
= \I t\, \phi(a_{t}), 
\\
& W \cdot \phi(a_{t}) = 0 \quad 
\mbox{for other basis vectors $W= X_{i}, Y_{i}, Z$ of $\lier{n}$}, 
\\
& W \cdot \phi(a_{t}) = -\tau_{\lambda^{\ast}}(W) \phi(a_{t}) \quad 
\mbox{for $W \in \lie{k}$}. 
\end{split}
\end{equation}
Here, we used the definition \eqref{eq:definition of psi} of
non-degenerate character $\eta$. 

Fix a non-degenerate invariant bilinear form $\langle \enskip, \enskip
\rangle$ on $\lier{g}$ and choose an orthonormal basis $\{W_{i}\}$ of
$\lier{s}$. 
Let $\pr_{\alpha^{\ast}}$ be the natural projection 
from 
$V_{\lambda^{\ast}}^{K} \otimes \lie{s} 
\simeq \oplus_{\alpha \in \Delta_{\lie{s}}} m(\alpha)\,
V_{(\lambda+\alpha)^{\ast}}^{K}$ to 
$V_{(\lambda+\alpha)^{\ast}}^{K}$. 
Then the $K$-type shift $\psi \mapsto \tilde{\psi}_{\alpha}$ which is
explained in \S\ref{subsection:shift operators} is translated into the
following operator: 
\begin{align*}
& 
P_{\alpha} : 
C^{\infty}(A \rightarrow V_{(\lambda/\gamma)^{\ast}}^{K/M^{\eta}})
\rightarrow 
C^{\infty}
(A \rightarrow V_{(\lambda+\alpha/\gamma)^{\ast}}^{K/M^{\eta}})
\\
& P_{\alpha} \phi(a_{t})
:=
\pr_{\alpha^{\ast}} \circ \nabla \phi(a_{t}), 
\qquad 
\nabla \phi(a_{t})
:= 
\sum_i W_i \cdot \phi(a_{t})\otimes W_i. 
\end{align*}
Actually, if we write 
$\iota_{\alpha}(v_{\alpha}) 
= \sum_{i} v_{\alpha}^{(i)} \otimes W_{i}$ 
for $v_{\alpha} \in V_{\lambda+\alpha}^{K}$, 
then 
\begin{align*}
\tilde{\psi}_{\alpha}(v_{\alpha})(a_{t}) 
&= \tilde{\psi}(\iota_{\alpha}(v_{\alpha}))(a_{t}) 
= \tilde{\psi}
(\sum_{i} v_{\alpha}^{(i)} \otimes W_{i})(a_{t}) 
= \sum_{i} 
L_{W_{i}} 
\psi(v_{\alpha}^{(i)})(a_{t}) 
\\
&= \sum_{i} 
\langle\langle v_{\alpha}^{(i)}, W_{i} \cdot \phi(a_{t})
\rangle\rangle_{\lambda} 
= 
\langle\langle 
\sum_{i} v_{\alpha}^{(i)} \otimes W_{i}, 
\sum_{j} W_{j} \cdot \phi(a_{t}) \otimes W_{j}
\rangle\rangle_{\lambda}' 
\\
&= 
\langle\langle 
\iota_{\alpha}(v_{\alpha}), 
\nabla \phi(a_{t}) 
\rangle\rangle_{\lambda}' 
\\
&= 
\langle\langle 
v_{\alpha}, 
P_{\alpha} \phi(a_{t})
\rangle\rangle_{\lambda+\alpha}. 
\end{align*}
Here, 
$\langle\langle Q' \otimes W_{i}, 
Q \otimes W_{j} \rangle\rangle_{\lambda}' 
:= \langle\langle Q', Q \rangle\rangle_{\lambda} 
\times \langle W_{i}, W_{j} \rangle$. 

In our $G=U(n,1)$ case, 
$\Delta_{\lie{s}} = \{\pm (\epsilon_{k}-\epsilon_{n+1})\, |\,
k=1,\dots,n\}$ and 
$(\lambda \pm (\epsilon_{n+1-k}-\epsilon_{n+1}))^{\ast} 
= \lambda^{\ast} \mp (\epsilon_{k}-\epsilon_{n+1})$. 
We write $P_{k}^{\pm}$ instead of 
$P_{\pm(\epsilon_{n+1-k}-\epsilon_{n+1})}$, for simplicity. 
These operators are calculated in \cite{T}. 
(The expression is slightly different because of the change of
notation and setting.)
\begin{proposition}\label{proposition:K-type shift}
Suppose $\phi(a_{t}) 
= \sum_{Q \in GT((\lambda/\gamma)^{\ast})} c(Q; t)\, Q 
\in 
C^\infty(A \rightarrow V_{(\lambda/\gamma)^{\ast}}^{K/M^{\eta}})$. 
Then the $K$-type shift operators $P_{k}^{\pm}$, $k=1,\dots,n$ 
is given by the following formulas: 
\begin{align}
P_{k}^{+} \phi(a_{t}) 
&= 
\sum_{Q \in GT((\lambda/\gamma)^{\ast})} 
b_{kn}(Q) 
(\theta-\vert \Lambda \vert - \vert \vect{q}_{n-1} \vert 
- 2 l_{k,n}-2n)\, c(Q; t)\, \sigma_{k,n}^{-}Q 
\label{eq:P_k^-}
\\
& \quad 
- 
t 
\sum_{i=1}^{n-1}
\sum_{\sigma_{i,n-1}^{+}Q \in GT((\lambda/\gamma)^{\ast})} 
\frac{b_{kn}(Q)\, a_{i,n-1}(Q)}{l_{k,n}-l_{i,n-1}} 
c(\sigma_{i,n-1}^{+}Q; t)\, \sigma_{k,n}^{-}Q, 
\notag
\\
P_{k}^{-} \phi(a_{t}) 
&= 
\sum_{Q \in GT((\lambda/\gamma)^{\ast})} 
a_{kn}(Q) 
(\theta+\vert \Lambda \vert + \vert \vect{q}_{n-1} \vert 
+ 2 l_{k,n}+2)\, c(Q; t)\, \sigma_{k,n}^{+}Q
\label{eq:P_k^+}
\\
& \quad 
+ 
t 
\sum_{i=1}^{n-1}
\sum_{\sigma_{i,n-1}^{-} Q \in GT((\lambda/\gamma)^{\ast})} 
\frac{a_{kn}(Q)\, b_{i,n-1}(Q)}{l_{k,n}-l_{i,n-1}+1} 
c(\sigma_{i,n-1}^{-}Q; t)\, \sigma_{k,n}^{+}Q. 
\notag
\end{align}
\end{proposition}

In order to state the next lemma, 
let $q_{0,n-2}:=\infty$ and $q_{n-1,n-2}:=-\infty$. 
\begin{lemma}\label{lemma:condition for vanishing under shift}
Let $\phi$ be an element of 
$C^\infty(A \rightarrow V_{(\lambda/\gamma)^{\ast}}^{K/M^{\eta}})$. 
\begin{enumerate}
\item
If $\gamma_{k-1}^{\ast} < \lambda_{k}^{\ast}$ 
and $P_{k}^{+} \phi = 0$ for $k \in \{2,\dots,n\}$, 
then $\phi = 0$. 
\item
If $\gamma_{k-1}^{\ast} > \lambda_{k}^{\ast}$ 
and $P_{k}^{-} \phi = 0$ for $k \in \{1,\dots,n-1\}$, 
then $\phi = 0$. 
\end{enumerate}
\end{lemma}
\begin{proof}
For any number $\ast$ depending on $Q \in GT(\lambda)$, 
we denote it by $\ast(Q)$, if we need to specify $Q$. 
For example, $q_{i,j}(Q)$ is the $q_{i,j}$ part of $Q \in
GT(\lambda)$. 

Since the proofs of these two are analogous, we shall show only (2). 
Let $Q_{0}$ be an element of $GT((\lambda/\gamma)^{\ast})$ which
satisfies $q_{k,n-1}(Q_{0}) = \lambda_{k}^{\ast}$, and let 
$Q_{1} := \sigma_{k,n-1}^{+}Q_{0}$. 
Then $Q_{1}$ is not in $GT((\lambda/\gamma)^{\ast})$, 
but $\sigma_{k,n-1}^{-}Q_{1}  =Q_{0} \in GT((\lambda/\gamma)^{\ast})$
and 
$\sigma_{k,n}^{+}Q_{1} \in GT(\lambda^{\ast}+(\epsilon_{k}-\epsilon_{n+1}))$, 
because $\gamma_{k-1}^{\ast} > \lambda_{k}^{\ast}$ implies that 
$\sigma_{k,n}^{+}Q_{1}$ satisfies the conditions in
  Definition~\ref{definition:GT}~(3): 
$q_{k-1,n-1}(Q_{1}) \geq 
\gamma_{k-1}^{\ast} \geq \lambda_{k}^{\ast}+1 
= q_{k,n}(\sigma_{k,n}^{+}Q_{1}) 
= q_{k,n-1}(\sigma_{k,n}^{+}Q_{1})$. 
Therefore, the term $\sigma_{k,n}^{+}Q_{1}$ appears in \eqref{eq:P_k^+}, 
and its coefficient in \eqref{eq:P_k^+} is 
\[
\frac{a_{k,n}(Q_{1})\, b_{k,n-1}(Q_{1})}
{l_{k,n}(Q_{1})-l_{k,n-1}(Q_{1})+1} 
c(\sigma_{k,n-1}^{-}Q_{1}; t) 
= 
\frac{a_{k,n}(Q_{0})\, b_{k,n-1}(\sigma_{k,n}^{+} Q_{1})}
{l_{k,n}(Q_{1})-l_{k,n-1}(Q_{1})+2} 
c(Q_{0}; t). 
\]
Here, we used the definition \eqref{eq:coefficients of GT} 
of $a_{i,j}(Q)$ and $b_{i,j}(Q)$. 
Since the coefficient of the right hand side is not zero, 
$c(Q_{0}; t)=0$ when $P_{k}^{-} \phi = 0$. 
We have shown that $c(Q; t)$ is zero for those $Q$ such that 
$q_{k,n-1}(Q)=\lambda_{k}^{\ast}$. 
By induction on $\lambda_{k}^{\ast}-q_{k.n-1}(Q)$ 
and by using \eqref{eq:P_k^+}, 
we can show $c(Q; t) = 0$ for
all $Q$, if $\gamma_{k-1}^{\ast}>\lambda_{k}^{\ast}$. 
\end{proof}


\subsection{Central elements of $U(\lie{gl}_{n+1})$} 
\label{subsection:central element}
In order to show Lemma~\ref{lemma:central element coming from shift}
below, we use the explicit forms of the elements in $Z(\lie{g})$. 
One of the most useful forms of the central elements of 
$U(\lie{gl}_{n+1})$ is the determinant type one (\cite{Capelli}). 
For the standard generator $E_{ij}$ of $\lie{gl}_{n+1}$ and a
parameter $u \in \C$, 
let $E_{ij}(u) := E_{ij}+u \delta_{ij}$ 
(Kronecker's delta). 
We define 
\[
C_{n+1}(u) 
:= 
\sum_{\sigma \in \mathfrak{S}_{n+1}} 
\sgn(\sigma) E_{\sigma(n+1),n+1}(u+n) E_{\sigma(n),n}(u+n-1) 
\cdots E_{\sigma(1),1}(u). 
\]
Then $C_{n+1}(u)$ is an element of $Z(\lie{gl}_{n+1})$ for
any $u$, 
and we obtain all the generators of $Z(\lie{gl}_{n+1})$ by
specializing $u$.  
Since 
$C_{n+1}(u) 
\equiv 
\prod_{p=1}^{n+1}(E_{pp}+u+p-1)$ modulo the left ideal generated by
strictly lower triangular matrices, 
the infinitesimal character is 
$\chi_{\Lambda}(C_{n+1}(u)) 
= 
\prod_{p=1}^{n+1}(\Lambda_{p}+(\rho_{\lie{g}})_{p}+u+p-1) 
= 
\prod_{p=1}^{n+1}(u+\Lambda_{p}+n/2)$. 
Here, 
$\rho_{\lie{g}} = \sum_{p=1}^{n+1} (\rho_{\lie{g}})_{p} \epsilon_{p} 
:= \frac{1}{2} \sum_{p=1}^{n+1} (n+2-2p) \epsilon_{p}$. 
\begin{lemma}\label{lemma:inf char of C_n+1}
$C_{n+1}(u)$ acts on $\stWhLgmg$ by the scalar 
$\prod_{p=1}^{n+1}(u+\Lambda_{p}+n/2)$. 
\end{lemma}

The exterior calculus is very useful for the manipulation of
non-commutative determinants. 
We use the method developed in \cite{IU}. 

The exterior algebra $\wedge \C^{2(n+1)}$ is an associative algebra
generated by $2(n+1)$ elements 
$e_{1}, \dots, e_{n+1}, 
e_{1}', \dots, e_{n+1}'$ subject to the relations 
$e_{i} e_{j} + e_{j} e_{i} = 0$, 
$e_{i}' e_{j} + e_{j} e_{i}' = 0$ and 
$e_{i}' e_{j}' + e_{j}' e_{i}' = 0$. 
We will work in the algebra 
$\wedge\C^{2(n+1)} \otimes U(\lie{gl}_{n+1})$,
where the subalgebras $\wedge\C^{2(n+1)}$ and $U(\lie{gl}_{n+1})$
commute with each other. 
Consider the following elements: 
\begin{align*}
& \eta_{j}(u) 
= 
\sum_{p=1}^{n+1} e_{p} E_{pj}(u), 
\qquad 
\eta_{i}'(u) 
= 
\sum_{q=1}^{n+1} e_{q}' E_{iq}(u), 
\\
& 
\Xi(u) = \sum_{p,q=1}^{n+1} e_{p} e_{q}' E_{pq}(u) 
= \sum_{j=1}^{n+1} \eta_{j}(u) e_{j}' 
= \sum_{i=1}^{n+1} e_{i} \eta_{i}'(u). 
\end{align*}
\begin{lemma}[\cite{IU}]
\label{lemma:commutation in exterior algebra}
\begin{enumerate}
\item
\label{eq:commutation, eta-eta}
For $i,j=1,\dots,n+1$, 
\begin{align*}
&\eta_{i}(u+1) \eta_{j}(u) + \eta_{j}(u+1) \eta_{i}(u) = 0,
&
&\eta_{i}'(u) \eta_{j}'(u+1) + \eta_{j}'(u) \eta_{i}'(u+1) = 0,
\end{align*}
\item
\label{eq:commutation, E_ij-omega}
For $i,j=1,\dots,n+1$, 
$[E_{ij}, \Xi(u)] 
= 
e_{j} \eta_{i}'(u)-\eta_{j}(u) e_{i}'$. 
\item\label{lemma:commutativity of Xi} 
For any $u,v \in \C$, $\Xi(u)$ and $\Xi(v)$ are commutative. 
\end{enumerate}
\end{lemma}
For $k \geq 0$, we consider the element 
\[
\Xi^{(k)}(u) = 
\Xi(u) \Xi(u-1) \cdots \Xi(u-k+1) 
= 
\Xi(u-k+1) \Xi(u-k+2) \cdots \Xi(u). 
\]
By Lemma~\ref{lemma:commutation in exterior algebra}, 
it is not hard to see that 
\begin{align}
\Xi^{(k)}(u) 
= & k \eta_{j}(u) e_{j}' \Xi^{(k-1)}(u-1) 
\label{eq:expansion-1}\\
& + (\Xi(u) - \eta_{j}(u) e_{j}') \cdots 
(\Xi(u-k+1) - \eta_{j}(u-k+1) e_{j}'), 
\notag\\ 
= & k e_{i} \eta_{i}'(u-k+1)  \Xi^{(k-1)}(u) 
\label{eq:expansion-2}\\
& + (\Xi(u-k+1) - e_{i} \eta_{i}'(u-k+1)) \cdots 
(\Xi(u) - e_{i} \eta_{i}'(u)),  
\notag
\\
0 = & \eta_{i}(u+n) e_{j}' \Xi^{(n)}(u+n-1) 
\quad 
\mbox{if $i\not=j$.} 
\label{eq:expansion-3}
\\
0 = & e_{i} \eta_{j}'(u) e_{k} e_{n+1}' \Xi^{(n-1)}(u+n-1) 
\quad 
\mbox{if $j\not=i, k$.} 
\label{eq:expansion-4}
\end{align}
From \eqref{eq:expansion-2}, we obtain 
\begin{equation}\label{eq:second expansion}
e_{j} e_{n+1}' \Xi^{(n)}(u+n-1) 
= n e_{j} e_{i} \eta_{i}'(u) e_{n+1}' \Xi^{(n-1)}(u+n-1). 
\end{equation}
By Proposition~2.2 in \cite{IU}, 
\begin{equation}\label{eq:determinant type}
\Xi^{(n+1)}(u+n) 
= (n+1)! C_{n+1}(u)\, \wedge^{\mathrm{top}},  
\end{equation}
where 
$\wedge^{\mathrm{top}}
:=e_{1} e_{1}' e_{2} e_{2}' \cdots e_{n+1} e_{n+1}'$. 
We need the cofactor expansion of $C_{n+1}(u)$ along the $(n+1)$-st
row and column. 

By \eqref{eq:basis of n}, we have 
\begin{align} 
\eta_{n+1}(u+n) - \eta_{n}(u+n) 
= & \sum_{p=1}^{n-1} \frac{1}{2} (-X_{p}+\I Y_{p}) e_{p}
+\frac{\I}{2} Z (e_{n}+e_{n+1})
\label{eq:eta_n+1-eta_n}\\
& +\frac{1}{2}(h-E_{nn}-E_{n+1,n+1}-2u-2n)(e_{n}-e_{n+1}),  
\notag\\
\eta_{n+1}'(u) + \eta_{n}'(u) 
= & \sum_{p=1}^{n-1} \frac{1}{2} (-X_{p}-\I Y_{p}) e_{p}'
+\frac{\I}{2} Z (-e_{n}+e_{n+1}) 
\label{eq:eta_n+1'-eta_n'}\\
& +\frac{1}{2}(h+E_{nn}+E_{n+1,n+1}+2u)(e_{n}'+e_{n+1}'). 
\notag
\end{align}
For two elements $x,y \in U(\lie{g})$, 
$x \equiv y$ means that they are equivalent modulo the right ideal
generated by $X-\eta_{t}(X)$, $X \in \lier{n}$. 
Since the actions of elements in $\lier{g}$ on 
$C^{\infty}(A \rightarrow V_{(\lambda/\gamma)^{\ast}}^{K/M^{\eta}})$ 
are given by \eqref{eq:action of Lie algebra elements}, 
we have 
\begin{align*}
& \frac{4}{(n+1) n} \Xi^{(n+1)}(u+n)
\notag\\
&= 
\frac{4}{n} (\eta_{n+1}(u+n)-\eta_{n}(u+n)) e_{n+1}' \Xi^{(n)}(u+n-1)
\\
&\equiv 
\frac{2}{n} 
\{-t e_{n-1} + (h-E_{nn}-E_{n+1,n+1}-2u-2n)(e_{n}-e_{n+1})\} 
e_{n+1}' \Xi^{(n)}(u+n-1)
\\
&= 
2\{t e_{n-1} - (h-E_{nn}-E_{n+1,n+1}-2u-2n) e_{n}\} 
\\
& \qquad \quad \times 
(\eta_{n+1}'(u)+\eta_{n}'(u)) e_{n+1} e_{n+1}' \Xi^{(n-1)}(u+n-1) 
\\
&\equiv 
\{t^{2} e_{n-1} e_{n-1}' 
+ t(h+E_{nn}+E_{n+1,n+1}+2u) e_{n-1} e_{n}' 
\\
& \qquad 
- t(h-E_{nn}-E_{n+1,n+1}-2u-2n) e_{n} e_{n-1}' 
\\
& \qquad 
- (h-E_{nn}-E_{n+1,n+1}-2u-2n) (h+E_{nn}+E_{n+1,n+1}+2u) 
e_{n} e_{n}' 
\} 
\\
& \qquad \quad \times 
e_{n+1} e_{n+1}' \Xi^{(n-1)}(u+n-1). 
\end{align*}
Here, we used \eqref{eq:expansion-1}, \eqref{eq:expansion-3} for the
first equality, 
\eqref{eq:eta_n+1-eta_n} for the second equivalence, 
\eqref{eq:expansion-4}, \eqref{eq:second expansion} for the third
equality, 
and \eqref{eq:eta_n+1'-eta_n'} for the last equivalence. 

By \eqref{eq:determinant type}, we have 
\begin{align*}
& e_{n} e_{n}' e_{n+1} e_{n+1}' \Xi^{(n-1)}(u+n-1) 
= (n-1)!\, C_{n-1}(u+1)\, \wedge^{\mathrm{top}}. 
\end{align*}
From this equation, we see that  
\begin{align*}
& 
e_{n} e_{n-1}' e_{n+1} e_{n+1}' \Xi^{(n-1)}(u+n-1) 
= 
(n-1)!\, \ad(E_{n-1,n}) C_{n-1}(u+1)\, \wedge^{\mathrm{top}}, 
\\
& 
e_{n-1} e_{n}' e_{n+1} e_{n+1}' \Xi^{(n-1)}(u+n-1) 
= 
-(n-1)! \, \ad(E_{n,n-1}) C_{n-1}(u+1)\, \wedge^{\mathrm{top}}, 
\\
& 
e_{n-1} e_{n-1}' e_{n+1} e_{n+1}' \Xi^{(n-1)}(u+n-1) 
\\
& \qquad 
= 
(n-1)! \, 
\{1-\ad(E_{n-1,n}) \ad(E_{n,n-1})\} C_{n-1}(u+1)\, \wedge^{\mathrm{top}}. 
\end{align*}

We shall write down the action of $C_{n+1}(u)$ on 
$\phi(a_{t}) 
= 
\sum_{Q \in GT((\lambda/\gamma)^{\ast})} 
c(Q; t) Q 
\in C^{\infty}(A \rightarrow 
\Hom_{M^{\eta}}(V_{\lambda}^{K},V_{\gamma}^{M^{\eta}}))$. 

Let $\bar{\lie{u}}_{\lie{m}}$ be the nilpotent subalgebra of $\lie{m}$
consisting of upper triangular matrices. 
Since the $\vect{q}_{n-1}$ part of 
$Q 
= (\vect{q}_{n-2}, \vect{q}_{n-1}, \vect{q}_{n}) 
\in GT((\lambda/\gamma)^{\ast})$ is a highest weight of 
$V_{\lambda^{\ast}}^{K} |_{U(n-1)}$, 
an element $z \in Z(\lie{m})$ acts on $Q$ by 
$\tau_{\lambda^{\ast}}({}^{t}z) Q 
= \gamma_{2}'(z)(-\vect{q}_{n-1}) Q$. 
($\gamma_{2}'$ is the non-shifted Harish-Chandra map 
defined in the proof of Theorem~\ref{theorem:generic case}.) 
Here, we used \eqref{eq:action of Lie algebra elements} and 
defined 
${}^{t}(Z_{1} \cdots Z_{l})
:= (-Z_{l})\cdots(-Z_{1})$ for $Z_{p} \in \lie{g}$. 
Since 
$\gamma_{2}'(C_{n-1}(u+1)) 
= E_{n-1,n-1}(u+n-1) \cdots E_{11}(u+1)$, 
it acts on $Q$ by the scalar 
$\prod_{p=1}^{n-1} (u+p-q_{p,n-1}) = \prod_{p=1}^{n-1}(u-l_{p,n-1})
=: S(Q)$. 

For $z \in Z(\lie{g})$, let
$\mathcal{D}_{\lambda,\gamma}(z)$ be the differential operator on 
$C^\infty(A \rightarrow V_{(\lambda/\gamma)^{\ast}}^{K/M^{\eta}})$
defined by $z \cdot \phi = \mathcal{D}_{\lambda,\gamma}(z) \phi$. 
Bringing the above results together, we get the following formula.  
\begin{proposition}\label{proposition:action of C_n+1}
The action of $C_{n+1}(u)$ on 
$\phi(a_{t}) 
\in 
C^{\infty}(A \rightarrow V_{(\lambda/\gamma)^{\ast}}^{K/M^{\eta}})$ is
expressed as follows: 
\begin{align}
-4 \mathcal{D}_{\lambda,\gamma}(C_{n+1}(u))\, &\phi(a_{t})
\notag\\
=
\sum_{Q \in GT((\lambda/\gamma)^{\ast})} 
&S(Q)
\Bigg[ 
\left\{
(\theta-n)^{2}
-(\vert \Lambda \vert + \vert \vect{q}_{n-1} \vert +2u+n)^{2}
-A(Q) t^{2} 
\right\} 
c(Q; t) 
\notag\\
& -t 
\sum_{p=1}^{n-1} 
\frac{a_{p,n-1}(Q)}{u-l_{p,n-1}} 
(\theta+\vert \Lambda \vert + \vert \vect{q}_{n-1} \vert 
+ 1 + 2u) c(\sigma_{p,n-1}^{+}Q; t) 
\label{eq:action of C_n+1}\\
& +t 
\sum_{p=1}^{n-1} 
\frac{b_{p,n-1}(Q)}{u-l_{p,n-1}} 
(\theta-\vert \Lambda \vert - \vert \vect{q}_{n-1} \vert 
+ 1 -2n- 2u) c(\sigma_{p,n-1}^{-}Q; t) 
\notag\\
&
-t^{2} 
\sum_{\genfrac{}{}{0pt}{}{p,r=1}{p\not=r}}^{n-1} 
\frac{b_{p,n-1}(Q)\, a_{r,n-1}(\sigma_{p,n-1}^{-}Q)}
{(u-l_{p,n-1}) (u-l_{r,n-1})} 
c(\sigma_{p,n-1}^{-} \sigma_{r,n-1}^{+}Q; t)
\Bigg] Q,
\notag\\
A(Q) := &
1-
\sum_{p=1}^{n-1} 
\frac{a_{p,n-1}(Q)^{2}-b_{p,n-1}(Q)^{2}}
{u-l_{p,n-1}}.
\notag
\end{align}
\end{proposition}

\begin{lemma}\label{lemma:central element coming from shift} 
Let $\tau_{\lambda}$ be a $K$-type of $\stWhLgmg$. 
On the space 
$C^{\infty}(A \rightarrow V_{(\lambda/\gamma)^{\ast}}^{K/M^{\eta}})$, 
the operators $P_{k}^{-} \circ P_{k}^{+}$ and 
$P_{k}^{+} \circ P_{k}^{-}$, $k \in \{1,\dots,n\}$, are
central, namely 
\begin{align}
P_{k}^{-} \circ P_{k}^{+}
&= \mathcal{D}_{\lambda,\gamma}(C_{n+1}(l_{k,n})), 
&
P_{k}^{+} \circ P_{k}^{-} 
&= \mathcal{D}_{\lambda,\gamma}(C_{n+1}(l_{k,n}+1)). 
\end{align}
\end{lemma}
\begin{proof}
By Propositions~\ref{proposition:K-type shift}, 
\ref{proposition:action of C_n+1}, 
we know that we may show the identity 
\[
1 
= 
\sum_{p=1}^{n-1} 
\frac{a_{p,n-1}(Q)^{2}}{l_{k,n}-l_{p,n-1}} 
- 
\sum_{p=1}^{n-1} 
\frac{b_{p,n-1}(Q)^{2}}{l_{k,n}-l_{p,n-1}+1}. 
\]
This identity is obtained by comparing the coefficient of 
$\sigma_{k,n}^{+} Q$ in the identity 
$\tau_{\lambda^{\ast}}(E_{n,n+1})\, Q 
= \tau_{\lambda^{\ast}}([E_{n,n-1}, [E_{n-1,n}, E_{n,n+1}]]) Q$. 
\end{proof}

\subsection{Determination of composition series}
\label{subsection:determination of composition series} 

In this subsection, we determine the composition series of
$\stWhLgmg$. 
When something concerning the irreducible modules
$\overline{\pi}_{a,b}$ is described, 
the statement concerning them is assumed to be excluded 
if there is not
such module $\overline{\pi}_{a,b}$, i.e. if $a+b > n+1$. 

\begin{lemma}\label{lemma:zero shifts}
Suppose that $\gamma$ is given by 
\eqref{eq:condition for sigma}. 
If a pair $V_{1}$ and $V_{2}$ satisfies one of the following
conditions, then there is no non-zero $\lie{g}$-action in
$\stWhLgmg$ which sends $V_{1}$ to $V_{2}$: 
\begin{enumerate}
\item
$V_{1} \simeq \overline{\pi}_{i,j}$, 
$V_{2} \simeq \overline{\pi}_{i+a,j},
  \overline{\pi}_{i, j+b}$, 
$a, b = \pm 1$. 
\item
$V_{1} \simeq \overline{\pi}_{i+a,j}$, 
$V_{2} \simeq \overline{\pi}_{i+a,j+b}$, 
$a, b = \pm 1$. 
\item
$V_{1} \simeq \overline{\pi}_{i,j+b}$, 
$V_{2} \simeq \overline{\pi}_{i+a,j+b}$, 
$a, b = \pm 1$. 
\end{enumerate}
(Any double signs are allowed.) 
\end{lemma}

\begin{proof}
If there is a $\lie{g}$-action sending an element of 
$\overline{\pi}_{a,b}$ to $\overline{\pi}_{a',b'}$, 
then the $K$-spectra $\widehat{K}(\overline{\pi}_{a,b})$ and 
$\widehat{K}(\overline{\pi}_{a',b'})$ should be adjacent, i.e. there
should be $K$-types 
$\tau_{\lambda} \in \widehat{K}(\overline{\pi}_{a,b})$ 
and $\tau_{\lambda'} \in \widehat{K}(\overline{\pi}_{a',b'})$ such
that $\lambda-\lambda'$ is a weight of $\lie{s}$. 

Assume that $V_{1}$ and $V_{2}$ are isomorphic to
$\overline{\pi}_{i,j}$ and $\overline{\pi}_{i-1,j}$,
respectively, and $V_{1} \rightarrow V_{2}$.  
Two $K$-types $\lambda \in \widehat{K}(\overline{\pi}_{i,j})$ and 
$\lambda' \in \widehat{K}(\overline{\pi}_{i-1,j})$ are adjacent if and
only if $\lambda_{i} = \Lambda_{i}-n/2+i$,
$\lambda_{i}'=\Lambda_{i}-n/2+i-1$ 
and $\lambda_{p}=\lambda_{p}'$ for $p\not=i$. 
Recall the discussion in \S~\ref{subsection:explicit shift operators}. 
The $\lie{s}$-action which sends an element of
$V_{\lambda}^{K} \subset V_{1}$ 
to 
$V_{\lambda'}^{K} \subset V_{2}$ is realized by the
shift operator $P_{n+1-i}^{-} \phi$. 

Consider the shift $P_{n+1-i}^{+} \circ P_{n+1-i}^{-} \phi$. 
Lemma~\ref{lemma:central element coming from shift} asserts 
that this is equal to 
$\mathcal{D}_{\lambda,\gamma}
(C_{n+1}(l_{n+1-i,n}+1)) 
\phi 
= \chi_{\Lambda}(C_{n+1}(l_{n+1-i,n}+1)) \phi$. 
By Lemma~\ref{lemma:inf char of C_n+1}, 
we know that $\chi_{\Lambda}(C_{n+1}(l_{n+1-i,n}+1)) = 0$, 
since 
$l_{n+1-i,n}
=
\lambda_{n+1-i}^{\ast}-(n+1-i) 
= -\lambda_{i}-(n+1-i) 
= -(\Lambda_{i}-n/2+i)-(n+1-i) 
= -\Lambda_{i}-n/2-1$. 
Therefore, $P_{n+1-i}^{+} P_{n+1-i}^{-} \phi=0$. 
Now, 
$\lambda_{i}' = \Lambda_{i}-n/2+i-1 < \gamma_{i-1}$
implies $(\lambda')_{n+1-i}^{\ast} > \gamma_{n-i}^{\ast}$, 
and the condition $1 \leq i \leq n-1$ of 
Proposition~\ref{proposition:condition for sigma, integral} 
is paraphrased as $2 \leq n+1-i \leq n$. 
Under these conditions, $P_{n+1-i}^{+}$ is injective
by Lemma~\ref{lemma:condition for vanishing under shift}, 
so $P_{n+1-i}^{-} \phi=0$. 
But this contradicts $V_{1} \rightarrow V_{2}$. 
Therefore, $V_{1} \not\rightarrow V_{2}$. 
Other cases can be shown analogously. 
\end{proof}

\begin{corollary}\label{corollary:pi_ij is multiplicity free}
The multiplicity of $\overline{\pi}_{i,j}$ in 
$\stWhLgmg$ is one. 
\end{corollary}
\begin{proof}
We know that the socle of $\stWhLgmg$ is
isomorphic to $\overline{\pi}_{i,j}$ 
(Proposition~\ref{proposition:composition factors-1}). 
Assume that there exists a composition factor $V_{1}$ 
which is isomorphic to $\overline{\pi}_{i,j}$ but is not in the
socle. 
By Proposition~\ref{proposition:composition factors-1}, 
$\overline{\pi}_{a,b}$ is a composition factor of 
$\stWhLgmg$ only if $a=i, i\pm1$ and
$b=j, j\pm1$. 
By Theorem~\ref{theorem:K-spectra}, it is adjacent to
$V_{1}$ if and only if $|a-a'|+|b-b'|=1$. 
Therefore, there exits a composition factor $V_{2}$ which is
isomorphic to one of
$\overline{\pi}_{i\pm 1,j}$, $\overline{\pi}_{i,j\pm 1}$ 
such that $V_{1} \rightarrow V_{2}$. 
But we have shown in Lemma~\ref{lemma:zero shifts} 
that this is impossible. 
\end{proof}

\begin{lemma}\label{lemma:second floor}
Suppose that $\gamma$ is given by \eqref{eq:condition for sigma}. 
\begin{enumerate}
\item
The socle of $\stWhLgmg/\overline{\pi}_{i,j}$ is 
$\overline{\pi}_{i-1,j} \oplus \overline{\pi}_{i,j+1} 
\oplus \overline{\pi}_{i,j-1} \oplus \overline{\pi}_{i+1,j}$. 
\item
The multiplicities of 
$\overline{\pi}_{i\pm 1,j}, \overline{\pi}_{i,j\pm 1}$ 
in $\stWhLgmg$ are all one. 
\end{enumerate}
\end{lemma}
\begin{proof}
(1) 
As we stated in the proof of the previous corollary, 
a composition factor $V$ is adjacent to $\overline{\pi}_{i,j}$
only if $V$ is isomorphic to one of 
$\overline{\pi}_{i\pm 1,j}$, $\overline{\pi}_{i,j\pm 1}$. 
So only 
$\overline{\pi}_{i \pm 1,j}$, $\overline{\pi}_{i,j \pm 1}$ can
be a simple submodule of $\stWhLgmg/\overline{\pi}_{i,j}$. 
We know from Proposition~\ref{proposition:at least one} 
that the multiplicity of each of them 
in $\stWhLgmg$ is at least one. 

Choose a composition factor $V$ isomorphic to, say,
$\overline{\pi}_{i-1,j}$. 
Other cases can be shown analogously. 
Recall the proof of Lemma~\ref{lemma:zero shifts}. 
Let $\phi$ be the function which 
characterizes the non-zero $K$-type $\lambda'$ of 
$V \simeq \overline{\pi}_{i-1,j}$. 
Assume that there is no non-zero $\lie{s}$-action from $V$
to the unique simple submodule which is isomorphic to
$\overline{\pi}_{i,j}$. 
Since the multiplicity of $\overline{\pi}_{i,j}$ in
$\stWhLgmg$ is one, $P_{n+1-i}^{+} \phi = 0$. 
But we have seen in the proof of 
Lemma~\ref{lemma:zero shifts} that this implies $\phi=0$. 
This is a contradiction, 
so $V \rightarrow \overline{\pi}_{i,j}$. 

Assume that there are two composition factors 
$V_{1}, V_{2}$ in the socle of 
$\stWhLgmg/\overline{\pi}_{i,j}$, both of which are isomorphic to 
$\overline{\pi}_{i-1,j}$. 
Let $\phi_{k}$, $k=1,2$, be the functions which 
characterize the non-zero $K$-type $\lambda'$ of 
$V_{k}$, respectively. 
Then both of $P_{n+1-i}^{+} \phi_{k}$ characterize the same
$K$-type $\lambda$ of the unique simple submodule, and the
multiplicity of this $K$-type is one, there are constants
$c_{k}$ such that 
$c_{1} P_{n+1-i}^{+} \phi_{1} 
= c_{2} P_{n+1-i}^{+} \phi_{2}$. 
Since $P_{n+1-i}^{+}$ is injective in this case, 
this implies that $c_{1} \phi_{1} + c_{2} \phi_{2} = 0$. 
This contradicts the fact that $\phi_{k}$, $k=1,2$,
characterize the $K$-types $\lambda'$ of different factors
$V_{k}$, $k=1,2$. 
Therefore, the multiplicity of $\overline{\pi}_{i-1,j}$ in
the socle of $\stWhLgmg/\overline{\pi}_{i,j}$
is one. 

(2) Assume that there is a composition factor $V_{1}$
isomorphic to, say, $\overline{\pi}_{i-1,j}$ but not in the
second floor of $\stWhLgmg$. 
The irreducible modules which are adjacent to 
$\overline{\pi}_{i-1,j}$ are $\overline{\pi}_{i,j}$ and 
$\overline{\pi}_{i-1,j\pm 1}$. 
Therefore, there exists a composition factor $V_{2}$ in the
third or higher floor such that it is isomorphic to one of the
above and $V_{1} \rightarrow V_{2}$. 
But this is impossible since 
(i) the multiplicity of $\overline{\pi}_{i,j}$ is one and it
is located in the bottom, 
and 
(ii) 
$\overline{\pi}_{i-1,j} \not\rightarrow
\overline{\pi}_{i-1,j\pm 1}$ 
by Lemma~\ref{lemma:zero shifts}(2). 
\end{proof}

\begin{lemma}\label{lemma:third floor}
Suppose that $\gamma$ is given by \eqref{eq:condition for sigma}. 
\begin{enumerate}
\item
The socle of 
$(\stWhLgmg
/
\overline{\pi}_{i,j})
/
(\overline{\pi}_{i-1,j} \oplus \overline{\pi}_{i,j+1} 
\oplus \overline{\pi}_{i,j-1} \oplus \overline{\pi}_{i+1,j}
/
\overline{\pi}_{i,j})$ is 
$\overline{\pi}_{i-1,j-1} 
\oplus \overline{\pi}_{i-1,j+1} 
\oplus \overline{\pi}_{i+1,j-1} 
\oplus \overline{\pi}_{i+1,j+1}$. 
Moreover, the non-zero $\lie{s}$-actions from the third
floor to the second are 
$\overline{\pi}_{a, j\pm 1} 
\rightarrow 
\overline{\pi}_{a, j}$ 
and 
$\overline{\pi}_{i\pm 1, b} 
\rightarrow 
\overline{\pi}_{i, b}$, 
$a=i\pm 1$, $b=j\pm 1$. 
\item
The multiplicities of $\overline{\pi}_{i\pm 1,j\pm 1}$ in 
$\stWhLgmg$ are all one. 
\end{enumerate}
\end{lemma}

\begin{proof}
The proof is almost the same as that of 
Lemma~\ref{lemma:second floor}. 

Since 
(i) the multiplicities of 
$\overline{\pi}_{i,j\pm 1}, \overline{\pi}_{i\pm 1,j}$ 
and $\overline{\pi}_{i,j}$ are one, 
and 
(ii) they are in the first or second floor, 
the third floor is a direct sum of 
$\overline{\pi}_{i\pm 1,j\pm 1}$'s, 
and there is no higher floor in $\stWhLgmg$. 

For each $\overline{\pi}_{i\pm 1,j\pm 1}$, 
there exists at least one factor isomorphic to it 
in $\stWhLgmg$. 
Suppose, say, $V_{1} \simeq \overline{\pi}_{i-1,j+1}$ and 
$V_{2} \simeq \overline{\pi}_{i,j+1}$, 
the latter is in the second floor. 
Let $\phi$ be the function which characterizes a $K$-type of 
$V_{1}$ adjacent to $\overline{\pi}_{i,j+1}$. 
The shift operator sending $\phi$ to a $K$-type of 
$V_{2}$ is $P_{n+1-i}^{+}$. 
The proof of Lemma~\ref{lemma:zero shifts} says 
that this is injective. 
Therefore, $V_{1} \rightarrow V_{2}$. 
The uniqueness of the factor isomorphic to 
$\overline{\pi}_{i-1,j+1}$ is shown in the same way as in
the proof of the previous lemma. 
\end{proof}

We have obtained the following second main theorem of this paper. 
\begin{theorem}\label{theorem:main}
Suppose that $G=U(n,1)$ and the infinitesimal character $\Lambda$ is
regular integral. 
If the highest weight of 
$\sigma \in \widehat{M^{\eta}} 
\simeq \widehat{U(n-2)} \times \widehat{U(1)}$ satisfies 
\eqref{eq:condition for sigma} 
for some $i=1,\dots,n-1$, $j=2,\dots,n+1-i$, 
then the composition series of $\stWhLgmg$ is 
\begin{equation}\label{eq:composittion series of I}
\stWhLgmg 
\quad \simeq \quad 
\begin{xy}
(0,12)*{\overline{\pi}_{i-1,j+1}}="A_{11}",
(15,12)*{\overline{\pi}_{i-1,j-1}}="A_{12}",
(30,12)*{\overline{\pi}_{i+1,j+1}}="A_{13}",
(45,12)*{\overline{\pi}_{i+1,j-1}}="A_{14}",
(0,0)*{\overline{\pi}_{i-1,j}}="A_{21}",
(15,0)*{\overline{\pi}_{i,j+1}}="A_{22}",
(30,0)*{\overline{\pi}_{i,j-1}}="A_{23}",
(45,0)*{\overline{\pi}_{i+1,j}}="A_{24}",
(22.5,-12)*{\overline{\pi}_{i,j}}="A_{31}",
\ar "A_{11}";"A_{21}"
\ar "A_{11}";"A_{22}"
\ar "A_{12}";"A_{21}"
\ar "A_{12}";"A_{23}"
\ar "A_{13}";"A_{22}"
\ar "A_{13}";"A_{24}"
\ar "A_{14}";"A_{23}"
\ar "A_{14}";"A_{24}"
\ar "A_{21}";"A_{31}"
\ar "A_{22}";"A_{31}"
\ar "A_{23}";"A_{31}"
\ar "A_{24}";"A_{31}"
\end{xy}
\end{equation}
Here, if $i+j=n$ or $n+1$, the modules 
$\overline{\pi}_{a,b}$, $a+b>n+1$, are regarded to be zero and
the arrows starting from or ending at such modules are omitted. 
\end{theorem}

\section{Ending remark}
\label{section:ending remark} 
In this paper, we characterized the module $\stWhLsmg$ by the
conditions (1)--(3) in \S~\ref{section:introduction}. 
The condition (1) is imposed to make the modules $\WhL$, 
$\stWhLs$ and $\stWhLsmg$ suitably small,
i.e. $K$-admissible. 
The author thinks that it is interesting to investigate the
structure of modules which are characterized by other
conditions that make the modules in question $K$-admissible. 
For example, if the real rank of $G$ is one, 
we may replace the condition 
``$f$ is a joint eigenfunction of $Z(\lie{g})$'' in (1) 
with 
``$f$ is an eigenfunction of the Casimir operator and admits 
a generalized infinitesimal character''. 
Under the latter condition, the module is still $K$-admissible and
has finite length. 
For $G=U(n,1)$, the composition series of such module 
(with the trivial generalized infinitesimal character) is  
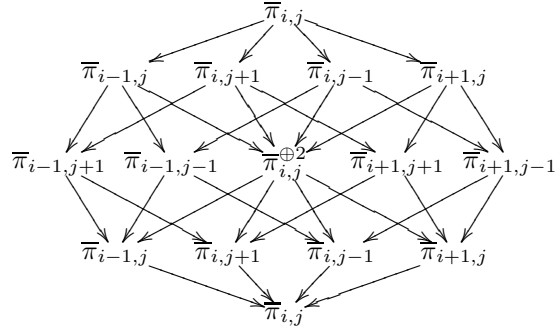
\[
\begin{xy}
(0,0)*{\overline{\pi}_{i-1,j+1}}="A_{11}",
(15,0)*{\overline{\pi}_{i-1,j-1}}="A_{12}",
(30,0)*{\overline{\pi}_{i,j}^{\oplus 2}}="A_{13}",
(45,0)*{\overline{\pi}_{i+1,j+1}}="A_{14}",
(60,0)*{\overline{\pi}_{i+1,j-1}}="A_{15}",
(7.5,-12)*{\overline{\pi}_{i-1,j}}="A_{21}",
(22.5,-12)*{\overline{\pi}_{i,j+1}}="A_{22}",
(37.5,-12)*{\overline{\pi}_{i,j-1}}="A_{23}",
(52.5,-12)*{\overline{\pi}_{i+1,j}}="A_{24}",
(30,-20)*{\overline{\pi}_{i,j}}="A_{31}",
(7.5,12)*{\overline{\pi}_{i-1,j}}="B_{21}",
(22.5,12)*{\overline{\pi}_{i,j+1}}="B_{22}",
(37.5,12)*{\overline{\pi}_{i,j-1}}="B_{23}",
(52.5,12)*{\overline{\pi}_{i+1,j}}="B_{24}",
(30,20)*{\overline{\pi}_{i,j}}="B_{31}",
\ar "A_{11}";"A_{21}"
\ar "A_{11}";"A_{22}"
\ar "A_{12}";"A_{21}"
\ar "A_{12}";"A_{23}"
\ar "A_{14}";"A_{22}"
\ar "A_{14}";"A_{24}"
\ar "A_{15}";"A_{23}"
\ar "A_{15}";"A_{24}"
\ar "A_{21}";"A_{31}"
\ar "A_{22}";"A_{31}"
\ar "A_{23}";"A_{31}"
\ar "A_{24}";"A_{31}"
\ar "B_{31}";"B_{21}"
\ar "B_{31}";"B_{22}"
\ar "B_{31}";"B_{23}"
\ar "B_{31}";"B_{24}"
\ar "B_{21}";"A_{11}"
\ar "B_{21}";"A_{12}"
\ar "B_{22}";"A_{11}"
\ar "B_{22}";"A_{14}"
\ar "B_{23}";"A_{12}"
\ar "B_{23}";"A_{15}"
\ar "B_{24}";"A_{14}"
\ar "B_{24}";"A_{15}"
\ar "B_{21}";"A_{13}"
\ar "B_{22}";"A_{13}"
\ar "B_{23}";"A_{13}"
\ar "B_{24}";"A_{13}"
\ar "A_{13}";"A_{21}"
\ar "A_{13}";"A_{22}"
\ar "A_{13}";"A_{23}"
\ar "A_{13}";"A_{24}"
\end{xy}
\]
This structure is more symmetric than that of $\stWhLsmg$.

\end{document}